\documentclass[a4paper,12pt]{amsart}
\usepackage{amsmath,amsfonts,amssymb,amsthm}
\usepackage{graphics}
\usepackage{epsfig}
\usepackage{esint}
\newcommand{\ignore}[1]{}

\newtheorem{proposition}{Proposition}

\newtheorem{lemma}{Lemma}
\newtheorem{corollary}{Corollary}

\newcommand{\R}{\mathbb{R}}

\usepackage{color}
\definecolor{darkred}{rgb}{0.9,0.1,0.1}
\definecolor{darkblue}{rgb}{0,0,0.7}
\definecolor{darkgreen}{rgb}{0,0.5,0}

\begin{document}
\begin{abstract}
These lecture notes present the
quantitative harmonic approximation result for quadratic optimal transport 
and general measures obtained in \cite{Goldman2021}. 
The aim is to give a clear presentation of the proof of \cite[Theorem 4.1]{Goldman2021} 
with more motivations, less PDE machinery, and a number of simplifications. 
\end{abstract}

\title[On the harmonic approximation to quadratic OT]{Lecture notes on the harmonic approximation to quadratic optimal transport}

\author{Lukas Koch, Felix Otto}

\maketitle




\section{A brief introduction}

These notes grew out of a couple of lecture series given by the second author
at 2022 summer schools on the 
topic of optimal transportation, its regularity theory, and its application to
the matching of random point clouds, see for instance
https://kantorovich.org/event/2022-optimal-transport-summer-school/schedule/. 
They presented results in
\cite{Goldman2020}, \cite{Goldman2021}, and \cite{Huesmann2021}. 

\medskip

The traditional approach to regularity \cite{Caffarelli1992} 
and partial regularity \cite{Figalli2010,DePhilippis2015} for optimal transportation relies
on the seminal regularity theory \cite{Caffarelli1990} 
for the corresponding Euler-Lagrange equation, the Monge-Amp\`ere equation,
based on the comparison principle.
The variational approach introduced in \cite{Goldman2020} avoids these arguments
and was first used to re-derive the partial regularity result of \cite{Figalli2010},
and then in \cite{Otto2021} to re-derive \cite{DePhilippis2015} for more general cost
functions of quadratic behavior.

\medskip

The added value of the variational approach lies in its robustness,
in particular in its ability to deal with general measures:
There is no need to have a Lebesgue density bounded away from zero and infinity in order
to allow for explicit barrier functions as in the approach based on comparison
principle. This robustness for instance allows to give a mesoscopic
characterization of the optimal transport between what is allowed to be an
atomic measure and the uniform distribution
\cite[Corollary 1.1]{Goldman2021}, and allows to analyze the matching between
independent copies of the Poisson point process \cite{Huesmann2021}.
We refer to www.mis.mpg.de/services/media/imprs-ringvorlesung-2022.html
for a gentle introduction into this aspect of matching.

\medskip

In analogy to de Giorgi's strategy for $\epsilon$-regularity of minimal surfaces,
the core of the variational regularity theory is a harmonic approximation result
\cite[Theorem 1.4]{Goldman2021},
which we will focus on in these notes.
Hence we do not discuss the literature further, 
but refer to \cite{Otto2021} for a careful review of the literature,
and the connection to minimal surface theory.
Compared to \cite[Section 3]{Goldman2021}, which we mainly rely on, 
these notes come with more motivations, less PDE machinery, and a couple of simplifications.
They allow for an independent reading.

\medskip

\subsection{Standing assumptions and language}

Throughout the entire text, we will consider two non-negative (finite) measures 
$\lambda$ and $\mu$ on $\R^d$ with $\lambda(\mathbb{R}^d)$ $=\mu(\mathbb{R}^d)$, 
which one should think of as two different spatial distributions of the same amount of mass. 
It is convenient to assume that $\lambda$ and $\mu$ have bounded support.
A non-negative measure $\pi$ on the product space $\mathbb{R}^d\times\mathbb{R}^d$
is called admissible if its marginals are given by $\lambda$ and $\mu$, 
which spelled out means
\begin{align}\label{ao00bis}
\int\zeta(x)\pi(dxdy)=\int\zeta d\lambda\quad\mbox{and}\quad
\int\zeta(y)\pi(dxdy)=\int\zeta d\mu
\end{align}
for all continuous and compactly supported functions (``test functions'') $\zeta$ on $\mathbb{R}^d$.
One should think of $\pi$ as one possible way of transporting the mass as distributed according 
to $\lambda$ into the shape as described by $\mu$ (a ''transport/transference plan'').
An admissible $\pi$ is called optimal if it minimizes 
\begin{align*}
\int|x-y|^2d\pi\stackrel{\mbox{short for}}{=}
\int_{\mathbb{R}^d\times\mathbb{R}^d}|x-y|^2\pi(dxdy),
\end{align*}
which one interprets as a total transportation cost, since it integrates the
cost of transporting a unit mass from $x$ to $y$, which here is given by the
{\it square} of the Euclidean distance. This is Kantorovich' relaxation of Monge's problem,
applied to the quadratic cost function. The infimum (which actually is attained)
\begin{align}\label{as01}
W^2(\lambda,\mu)=\inf\{\,\int|x-y|^2d\pi\,|\,\pi\;\mbox{is admissible for}\;\lambda,\mu\,\}
\end{align}
defines a distance function $W$, which we call Wasserstein distance.


\section{Connection of optimal transportation 
and the Neumann problem for the Poisson equation}\label{sec:orth}

In this section, we motivate the connection between optimal transportation (OT)
and the Neumann boundary value problem for the Poisson equation.

\subsection{Trajectories} 
For the above connection, it is convenient to adopt a dynamical view upon OT,
identifying a pair $(x,y)$ of (matched) points with the (straight) trajectory 
\begin{align}\label{ao00}
[0,1]\ni t\mapsto X(t):=ty+(1-t)x.
\end{align}
Given an optimal transfer plan $\pi$ for $\lambda,\mu$, we ask the question on how
to choose a function $\phi$ in such a way that its gradient $\nabla\phi$
captures the velocity of the trajectories, meaning
\begin{align}\label{ao65}
\dot X(t)\approx\nabla\phi(X(t))\quad\mbox{for}\;(x,y)\in{\rm supp}\pi.
\end{align}
As we shall see, the answer relates to the Poisson equation $-\triangle\phi$ $=\mu-\lambda$.

\medskip

\begin{figure}[h!]
  \centering
        \includegraphics[width=0.5\linewidth]{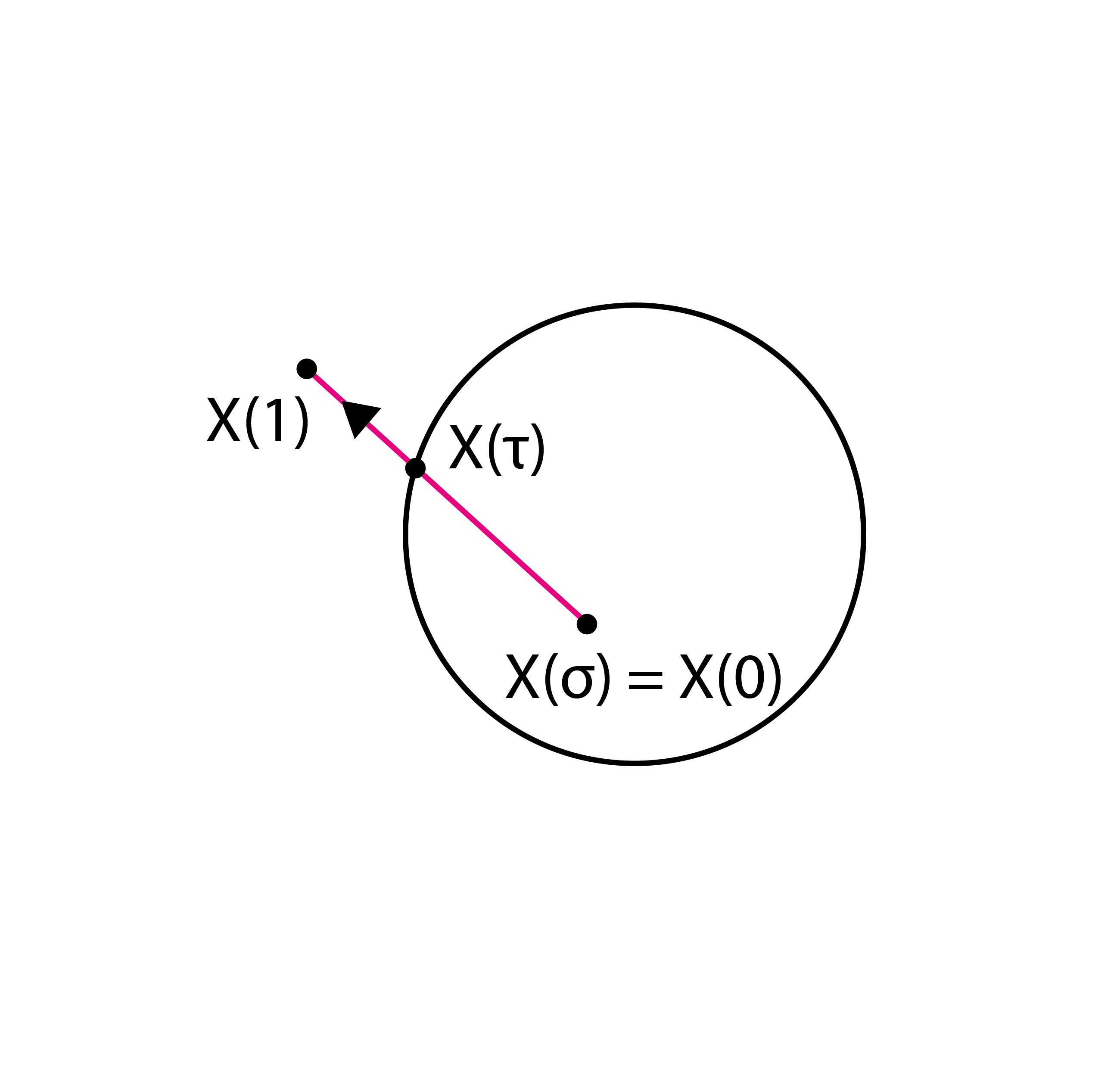}
    \caption{Entering and exiting times of trajectories}
    \label{fig:enteringExiting}
\end{figure}

We are interested in connecting to a boundary value problem 
for the Poisson equation on some domain,
say a ball $B_R$ of some radius $R$ (to be optimized later) 
and center w.~l.~o.~g.~given by the origin.
We are thus led to restrict ourselves\footnote{we will proceed to a further
restriction in (\ref{ao17})} to the set of trajectories
that spend some time in the closure $\bar B_R$:
\begin{align}\label{cw00}
\Omega:=\{\,(x,y)\,|\,\exists t\in[0,1]\;X(t)\in \bar B_R\,\}.
\end{align}
To every 
$(x,y)\in\Omega$, we associate the entering and exiting times $0\le\sigma\le\tau\le 1$
of the corresponding trajectory
\begin{align}\label{ao66}
\begin{array}{cc}
&\sigma:=\min\{t\in[0,1]\, |\, X(t)\in \bar B_R\},\\[1ex]
&\tau:=\max\{t\in[0,1]\,|\, X(t)\in \bar B_R\},
\end{array}
\end{align}
see also Fig. \ref{fig:enteringExiting}.
(Note that some trajectories may both enter and exit.)

Given a transfer plan $\pi$, 
we keep track of {\it where} the trajectories enter and exit $B_R$, which is captured by
two (non-negative) measures $f$ and $g$ concentrated on $\partial B_R$, defined through
\begin{align}
\int\zeta df&=\int_{\Omega\cap\{X(\sigma)\in\partial B_R\}} \zeta(X(\sigma))d\pi,\label{ao09}\\
\int\zeta dg&=\int_{\Omega\cap\{X(\tau)\in\partial B_R\}} \zeta(X(\tau))d\pi\label{ao09bis}
\end{align}
for all 
test functions functions $\zeta$. Note that the set of trajectories
$\Omega\cap\{X(\sigma)\in\partial B_R\}$
implicitly defines a Borel measurable 
subset of $\mathbb{R}^d\times\mathbb{R}^d$, namely the pre-image 
under the mapping (\ref{ao00}), which is continuous from $\mathbb{R}^d\times\mathbb{R}^d$
into $C^0([0,1])$.
Hence the integration against $\pi$ in (\ref{ao09}) is legitimate.

\begin{lemma}\label{lem:orth} 
We have for any admissible $\pi$ and any continuously differentiable function $\phi$ on $\bar B_R$
\begin{align}
\lefteqn{\int_\Omega\int_{\sigma}^\tau|\dot X(t)-\nabla\phi(X(t))|^2 dt d\pi}\nonumber\\
&=\int_\Omega\int_{\sigma}^\tau|\dot X(t)|^2 dt d\pi
+\int_\Omega\int_{\sigma}^\tau|\nabla\phi(X(t))|^2 dtd\pi\nonumber\\
&-2\int_{B_R}\phi d(\mu-\lambda)
-2\int_{\partial B_R} \phi d(g-f)\label{ao08}.
\end{align}
For later purpose, we record
\begin{align}\label{ao44}
\lambda(B_R)+f(\partial B_R)=\mu(B_R)+g(\partial B_R).
\end{align}
\end{lemma}


{\sc Proof of Lemma \ref{lem:orth}}.
For identity (\ref{ao08}) we note that for the mixed term we have by the chain rule
$\dot X(t)\cdot\nabla\phi(X(t))=\frac{d}{dt}[\phi(X(t))]$ and thus by the fundamental
theorem of calculus
$\int_\sigma^\tau\dot X(t)\cdot\nabla\phi(X(t))dt$ $=\phi(X(\tau))$ $-\phi(X(\sigma))$.
In view of definition (\ref{ao66}) we either have $X(\sigma)\in\partial B_R$ or
$X(\sigma)\in B_R$. By definition (\ref{ao66}) of $\sigma$ the latter implies
$\sigma=0$ and thus $X(\sigma)=x$,
so that the constraint $(x,y)\in\Omega$ may be dropped.
Hence
$\int_{\Omega}\phi(X(\sigma)) d\pi$
$=\int_{\Omega\cap\{X(\sigma)\in\partial B_R\}}\phi(X(\sigma)) d\pi$
$+\int_{\{x\in B_R\}}\phi(x) d\pi$. By definition (\ref{ao09}),
the first integral is $\int\phi df$. By admissibility (\ref{ao00bis}) of $\pi$,
the second integral is $\int_{B_R}\phi d\lambda$.
Likewise, one obtains $\int_{\Omega}\phi(X(\tau)) d\pi$
$=\int\phi dg+\int_{B_R}\phi d\mu$.

\medskip

Specifying to $\phi=1$, and thus $\nabla\phi=0$ so that the mixed term vanishes,
we learn (\ref{ao44}) from the above two identities.


\subsection{Perturbative regime}\label{ss:pert}
We will focus on a ``perturbative regime'', which comes in form of two local smallness conditions.
Any smallness condition has to be formulated in a non-dimensionalized way, 
which we implement by expressing this local smallness condition on a ball of non-dimensionalized
radius, it will be convenient to take $5$ as this radius. 

\medskip

The first smallness condition involves
the data (thus the letter $D$), that is, the two measures $\lambda$ and $\mu$. We monitor how close 
these measures are to the Lebesgue measure on $B_5$. It is natural to quantify
this in terms of the Wasserstein distance, see (\ref{as01}). Since the mass $\lambda(B_5)$
in general is not equal to the Lebesgue volume $|B_5|$, we have to split this into two:
We monitor how Wasserstein-close the restriction $\lambda\llcorner B_5$ is to
the uniform measure $\kappa_{\lambda}dx\llcorner B_5$, 
where $\kappa_{\lambda}$ $:=\frac{\lambda(B_5)}{|B_5|}$, and we monitor how close this
density $\kappa_{\lambda}$ is to unity. It is convenient to do both on the squared level:
\begin{align}\label{ao88}
D:&=W^2(\lambda\llcorner B_5,\kappa_\lambda dx\llcorner B_5)+(\kappa_\lambda-1)^2\nonumber\\
&+\mbox{same expression with $\lambda\leadsto\mu$}.
\end{align}

\medskip

In view of the localization (\ref{ao88}), it is convenient to further restrict the set
of trajectories, imposing that they start or end in $B_4$, 
thereby replacing (\ref{cw00}) by 
\begin{align}\label{ao17}
\Omega = \{(x,y)\in (B_4\times \R^d)\cup (\R^d\times B_4)|\; \exists t\in[0,1]\; X(t)\in \bar B_R\}.
\end{align}

\medskip

The second smallness condition involves the solution itself, i.~e.~$\pi$. 
It monitors the length of
trajectories that start or end in $B_5$. It does so in a square-averaged sense,
like the total cost function itself. In fact, it is a localization of the cost
functional (or energy, thus the letter $E$):
\begin{align}\label{ao45}
E:=\int_{(B_5\times\mathbb{R}^d)\cup(\mathbb{R}^d\times B_5)}|x-y|^2d\pi.
\end{align}

\medskip

We expect (and shall rigorously argue in Subsection \ref{PDEestimates}) 
that in the perturbative regime $E+D\ll 1$\footnote{Here and in the following we use the notation $\ll 1$ to mean that there is $\varepsilon>0$ such that the statement holds if $(E+D)\leq \varepsilon$.} and for a suitable $R\in[2,3]$
we have for the second r.~h.~s.~term in (\ref{ao08})
\begin{align}\label{ao14}
\int_\Omega\int_{\sigma}^\tau|\nabla\phi(X(t))|^2 dtd\pi
\approx\int_{B_R}|\nabla\phi|^2.
\end{align}
Indeed, for $E\ll 1$, trajectories are short so that 
\begin{align*}
\int_\Omega\int_{\sigma}^\tau|\nabla\phi(X(t))|^2 dtd\pi
\approx\int_{\{x\in B_R\}}|\nabla\phi(x)|^2 d\pi=\int_{B_R}|\nabla\phi|^2 d\lambda,
\end{align*}
where the last identity follows from admissibility (\ref{ao00bis}). 
Furthermore, for $D\ll 1$, $\lambda$ is close to Lebesgue
so that
\begin{align*}
\int_{B_R}|\nabla\phi|^2 d\lambda\approx\int_{B_R}|\nabla\phi|^2.
\end{align*}

\medskip


\subsection{Connection to the Neumann problem for the Poisson equation}
Hence in order to achieve (\ref{ao65}), in view of (\ref{ao08})
and (\ref{ao14}), we are led to minimize
\begin{align}\label{ao11}
\int_{B_R}|\nabla\phi|^2
-2\int_{B_R}\phi d(\mu-\lambda)-2\int_{\partial B_R}\phi d(g-f)
\end{align}
in $\phi$. 
A minimizer $\phi$ of (\ref{ao11}), 
if it exists as a continuously differentiable function on $\bar B_R$, 
would be characterized by the Euler-Lagrange equation 
\begin{align}\label{ao62}
\int_{B_R}\nabla\zeta\cdot\nabla\phi
-\int_{B_R}\zeta d(\mu-\lambda)-\int_{\partial B_R}\zeta d(g-f)=0
\end{align}
for all continuously differentiable test functions $\zeta$ on $B_R$.
If $\phi$ even exists as a twice continuously differentiable function on $\bar B_R$,
we could appeal to the calculus identity 
$\nabla\zeta\cdot\nabla\phi$ $=\nabla\cdot(\zeta\nabla\phi)-\zeta\triangle\phi$
and the divergence theorem in form of $\int_{B_R}\nabla\cdot(\zeta\nabla\phi)$
$=\int_{\partial B_R}\zeta\nu\cdot\nabla\phi$, where $\nu(x)=\frac{x}{R}$ denotes the outer 
normal to $\partial B_R$ in a point $x$, to obtain the integration by parts formula
\begin{align}\label{ao68}
\int_{B_R}\nabla\zeta\cdot\nabla\phi
=\int_{B_R}\zeta(-\triangle\phi)
+\int_{\partial B_R}\zeta\nu\cdot\nabla\phi.
\end{align}
Hence (\ref{ao62}) can be reformulated and regrouped as
\begin{align}\label{ao63}
\int_{B_R}\zeta(-\triangle\phi-d(\mu-\lambda))+\int_{\partial B_R}\zeta(\nu\cdot\nabla\phi
-d(g-f))=0.
\end{align}
Considering first all test functions $\zeta$'s that vanish on $\partial B_R$,
we learn from (\ref{ao63}) that $-\triangle\phi=\mu-\lambda$ distributionally in $B_R$.
Since $\mu-\lambda$ is a bounded measure, the first term in (\ref{ao63}) thus
vanishes also for test functions that do not vanish on $\partial B_R$.
Hence the second term in (\ref{ao63}) vanishes individually, which means
$\nu\cdot\nabla\phi=g-f$ distributionally on $\partial B_R$.
Hence we end up with what is called the Poisson equation with Neumann boundary conditions
\begin{align}\label{ao12}
-\triangle\phi=\mu-\lambda\;\mbox{in}\;B_R,\quad
\nu\cdot\nabla\phi=g-f\;\mbox{on}\;\partial B_R.
\end{align}
This is a classical elliptic boundary value problem, which for sufficiently regular $\mu-\lambda$
and $g-f$ has a unique twice differentiable solution, provided (\ref{ao44}) holds, and
\begin{align}\label{ao61}
\int_{B_R}\phi=0
\end{align}
is imposed. This motivates the connection between optimal transportation and
the (short) Neumann-Poisson problem.

\medskip

However, for rough (like sum of Diracs) measures $\lambda,\mu$, and thus also rough measures $f,g$,
the solution $\phi$ of (\ref{ao12}), even if it exists for this linear problem, 
will be rough, too. In particular, (\ref{ao14}) may not be true; even worse,
both the l.~h.~s.~and the r.~h.~s.~might be infinite. 
Hence we shall approximate both $\mu-\lambda$ and $g-f$ by smooth functions
(in fact, we shall approximate $\mu-\lambda$ by a constant function).
The best way to organize the output of Lemma \ref{lem:orth} is given by

\begin{corollary}\label{cor:orth}
We have for any admissible $\pi$ and any twice continuously differentiable 
function $\phi$ on $\bar B_R$
\begin{align}
\lefteqn{\int_\Omega\int_{\sigma}^\tau|\dot X(t)-\nabla\phi(X(t))|^2 dt d\pi}\nonumber\\
&\le\int_\Omega|x-y|^2 d\pi-\int_{B_R}|\nabla\phi|^2\nonumber\\
&+2\int_{B_R}\phi(-\triangle\phi-d(\mu-\lambda))
+2\int_{\partial B_R}\phi (\nu\cdot\nabla\phi-d(g-f))\nonumber\\
&+\int_\Omega\int_\sigma^\tau|\nabla\phi(X(t))|^2 dtd\pi
-\int_{B_R}|\nabla\phi|^2.\label{ao10}
\end{align}
\end{corollary}

As we argued, see (\ref{ao14}), we expect the term in last line (\ref{ao10})
to be of higher order. The integrals in the second r.~h.~s.~line can be made small
by approximately solving (\ref{ao12}) -- and there will be a trade-off
between making the last line and the second line small.
However, the main open task is to argue, based on the optimality of $\pi$, 
that the difference in the first r.~h.~s.~line 
is small for an approximate solution of (\ref{ao12}). Hence we turn to this task before 
dealing with the second and third line in Subsection \ref{PDEestimates}.

\medskip

{\sc Proof of Corollary \ref{cor:orth}}.
The upgrade of identity (\ref{ao08}) to inequality (\ref{ao10}) relies on 
\begin{align}
\int_{\Omega}\int_\sigma^\tau|\dot X(t)|^2dtd\pi&\le\int_{\Omega}|x-y|^2d\pi,\label{ao46}\\
\int_{B_R}|\nabla\phi|^2&=\int_{B_R}\phi(-\triangle\phi)+\int_{\partial B_R}\phi\nu\cdot\nabla\phi,
\label{ao48}
\end{align}
after adding and subtracting $2\int_{B_R}|\nabla \phi|^2$.
Inequality (\ref{ao46}) follows from 
$\int_\sigma^\tau|\dot X(t)|^2dt$ $\le\int_0^1|\dot X(t)|^2dt$ $=|x-y|^2$. 
Identity (\ref{ao48}) follows from (\ref{ao68}) for $\zeta=\phi$.

\medskip


\subsection{Localizing optimality}
As mentioned after Corollary \ref{cor:orth}, the main open task is to estimate
the first r.~h.~s.~line of (\ref{ao10}). For this, we will (for the first time) use that $\pi$ is
optimal. In order to connect to the Neumann-Poisson problem on $B_R$, we need to 
leverage optimality in a localized way. Of course, it will in general not be true that
the cost of $\pi$ localized to $(B_R\times\mathbb{R}^d)\cup(\mathbb{R}^d\times B_R)$
is estimated by the transportation cost between the localized measures
$\lambda\llcorner B_R$ and $\mu\llcorner B_R$. However, this is almost true
if one adds the distribution of the entering points $f$, see (\ref{ao09}), and exiting points $g$, see (\ref{ao09bis}),
respectively:

\begin{lemma}\label{lem:opt} For $\pi$ optimal we have
\begin{align}\label{ao70}
\big(\int_\Omega|x-y|^2d\pi\big)^\frac{1}{2}
&\le W(\lambda\llcorner B_R+f,\mu\llcorner B_R+g)\nonumber\\
&+\big(2\int_{\Omega\cap{\{\exists t\in[0,1]\;X(t)\in\partial B_R\}}}|x-y|^2d\pi\big)^\frac{1}{2}.
\end{align}
\end{lemma}

Lemma \ref{lem:opt} controls the transportation cost coming from those
trajectories that spend some time in $\bar B_R$, which amounts to the l.~h.~s.~of
(\ref{ao70}) according to definition (\ref{ao17}), by an OT problem localized to $\bar B_R$
as described by the first r.~h.~s.~term. It does so up to the transportation cost
coming from those (fewer) trajectories that cross (or touch) the boundary $\partial B_R$, see the
second r.~h.~s.~term. We shall argue in Lemma \ref{lem:linfty}, cf.~(\ref{m05}), 
that this last term (without the square root) is $o(E)$ for a good choice of $R$.

\medskip

As its form suggests, (\ref{ao70}) has the structure of a triangle inequality.
In fact, its proof has similarities with the proof of the triangle inequality for $W$,
using a disintegration (or conditioning) argument, c.~f.~\cite[Section 5.1]{Santambrogio}.

\medskip

\begin{figure}[h!]
  \centering
    \includegraphics[width=0.5\linewidth]{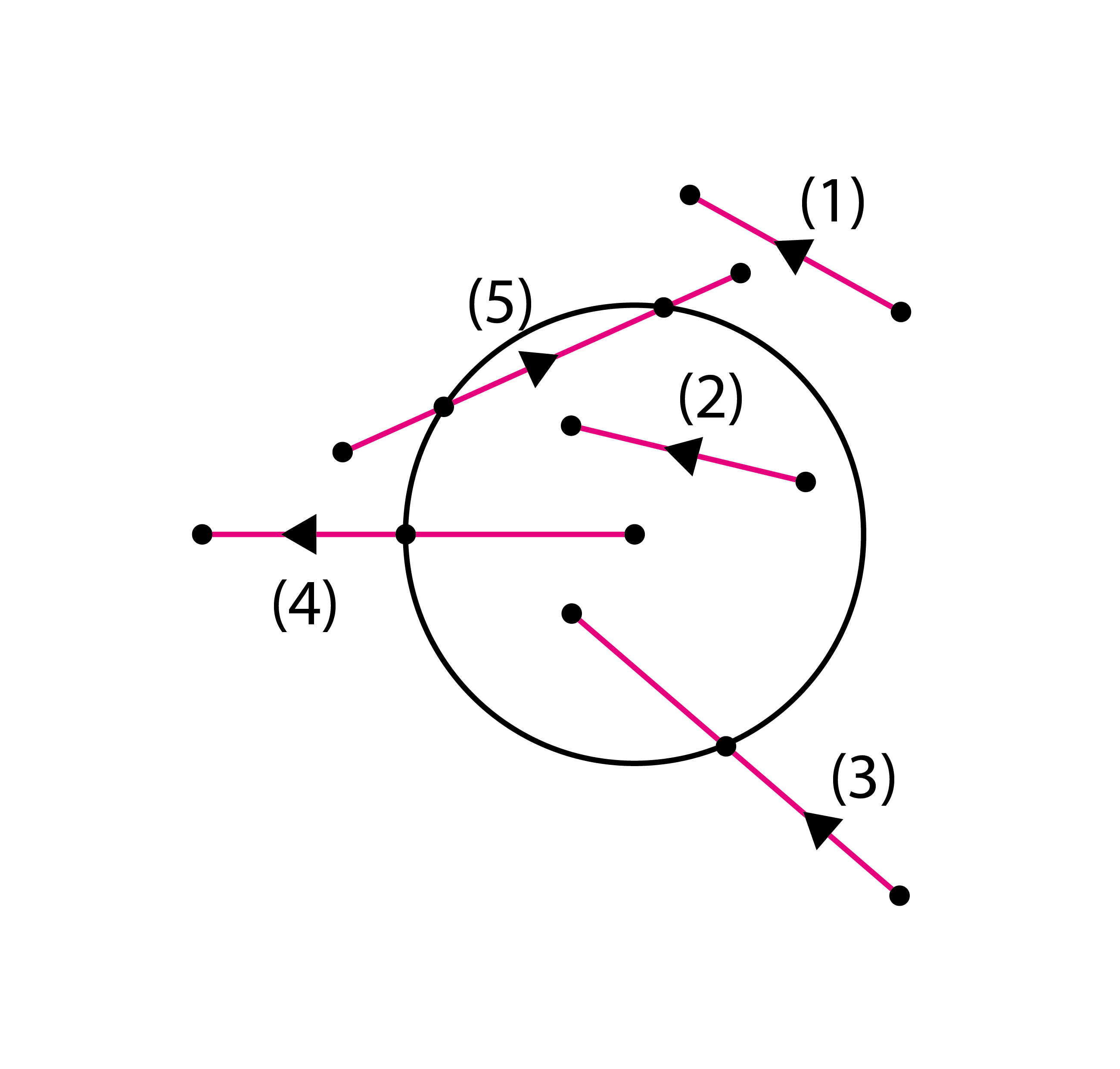}
    \caption{The terms in \eqref{ao73}}
\label{fig:terms}
\end{figure}

{\sc Proof of Lemma \ref{lem:opt}}. We now introduce the distribution of $x=X(0)$
under $\pi$ conditioned on the event that the trajectory $X$ enters at $z\in\partial B_R$.
In less probabilistic and more measure-theoretic terms (``dis-integration''), 
we introduce the (weakly continuous) family of probability measures 
$\{\lambda_z\}_{z\in\partial B_R}$ such that
\begin{align}\label{ao79}
\int_{\Omega\cap\{X(\sigma)\in \partial B_R\}}\zeta(x,X(\sigma))\pi(dxdy)
=\int_{\partial B_R}\int\zeta(x,z)\lambda_z(dx)f(dz),
\end{align}
which is possible by (\ref{ao09}). Here, $\zeta$ is an arbitrary test function on 
$\mathbb{R}^d\times\mathbb{R}^d$. Likewise, we introduce the probability
distribution $\{\mu_w\}_{w\in\partial B_R}$ of the end points of trajectories that
exit in $w$:
\begin{align}\label{ao80}
\int_{\Omega\cap\{X(\tau)\in \partial B_R\}}\zeta(X(\tau),y)\pi(dxdy)
=\int_{\partial B_R}\int\zeta(w,y)\mu_w(dy)g(dw).
\end{align}
Let $\bar\pi$ denote an optimal plan for $W(\lambda\llcorner B_R+f,\mu\llcorner B_R+g)$.
Equipped with these objects, we now define a competitor $\tilde\pi$ for $\pi$ 
that mixes $\pi$ with $\bar\pi$, in the sense that it takes the trajectories from
$\pi$ that stay outside of $\bar B_R$, the trajectories from $\bar\pi$ that stay
inside (the open) $B_R$, and concatenates trajectories $X$ from $\pi$ that enter or exit
$\bar B_R$ with trajectories of $\bar\pi$ that start or end in $\partial B_R$:
\begin{align}\label{ao73}
\int\zeta(x,y)\tilde\pi(dxdy)
&=\int_{\Omega^c}\zeta(x,y)\pi(dxdy)\nonumber\\
&+\int_{B_R\times B_R}\zeta(x,y)\bar\pi(dxdy)\nonumber\\
&+\int_{\partial B_R\times B_R}\int\zeta(x,y)\lambda_z(dx)\bar\pi(dzdy)\nonumber\\
&+\int_{B_R\times \partial B_R}\int\zeta(x,y)\mu_w(dy)\bar\pi(dxdw)\nonumber\\
&+\int_{\partial B_R\times \partial B_R}\int\int\zeta(x,y)\mu_w(dy)\lambda_z(dx)\bar\pi(dzdw)\nonumber\\
&=:(1) + (2) + (3) + (4) + (5),
\end{align}
see Fig. \ref{fig:terms}.

It is straightforward to see that $\tilde\pi$ has marginals $\lambda$ and $\mu$; 
by symmetry, it is sufficient to check
the first condition in (\ref{ao00bis})
by using (\ref{ao73}) for a function $\zeta=\zeta(x)$:
Since $\mu_w$ is a probability measure to the effect of $\int\zeta(x)\mu_w(dy)$ $=\zeta(x)$,
the second and the fourth r.~h.~s.~term of (\ref{ao73}) combine to
$\int_{B_R\times\mathbb{R}^d}\zeta(x)\bar\pi(dxdy)$ because $\bar\pi$,
like $\mu\llcorner B_R+g$, is supported in $\bar B_R$ (in $dy$). 
Likewise, the third and the fifth term combine to
$\int_{\partial B_R\times\mathbb{R}^d}\int\zeta(x)$ $\lambda_z(dx)$ $\bar\pi(dzdy)$.
By admissibility of $\bar\pi$, the combination of the second and fourth term
gives $\int_{B_R}\zeta(x)$ $\mu(dx)$, which as in the proof of Lemma \ref{lem:orth}
(by admissibility of $\pi$)
can be seen to be $\int_{\Omega\cap\{X(\sigma)\in B_R\}}\zeta(x)\pi(dxdy)$. 
Since $\int\zeta(x)\lambda_z(dx)$ does not depend on 
$y$, for the same reason, the combination of the third and fifth term
renders $\int_{\partial B_R}\zeta(z)f(dz)$, which by definition (\ref{ao09})
is equal to $\int_{\Omega\cap\{X(\sigma)\in \partial B_R\}}\zeta(x)$ $\pi(dxdy)$.
Hence these four terms combine to $\int_{\Omega}\zeta(x)\pi(dxdy)$.
Therefore, the r.~h.~s.~of (\ref{ao73}) collapses as desired to
$\int\zeta(x)\pi(dxdy)$, which coincides with $\int\zeta(x)\lambda(dx)$ by
admissibility of $\pi$.

\medskip

By optimality of $\pi$, we have $\int|x-y|^2d\pi\le\int|x-y|^2d\tilde\pi$;
rewriting this as $\int_\Omega|x-y|^2d\pi+\int_{\Omega^c}|x-y|^2d\pi\le\int|x-y|^2d\tilde\pi$,
and using (\ref{ao73}) for $\zeta(x,y)=|x-y|^2$, we gather
\begin{align}\label{ao81}
\big(\int_\Omega|x-y|^2d\pi\big)^\frac{1}{2}\le\|(f_2,f_3,f_4,f_5)\|,
\end{align}
where the four functions $f_2,\cdots,f_5\ge 0$ are given by
\begin{align*}
f_2(x,y):=|x-y|,\quad f_5^2(z,w):=\int\int|x-y|^2\mu_w(dy)\lambda_z(dx),\nonumber\\
f_3^2(z,y):=\int|x-y|^2\lambda_z(dx),\quad
f_4^2(x,w):=\int|x-y|^2\mu_w(dy),
\end{align*}
and the vector-valued $L^2$-type norm is defined through
%
\begin{align}\label{ao75}
\lefteqn{\|(f_2,f_3,f_4,f_5)\|^2}\nonumber\\
&=\int_{B_R\times B_R}f_2^2(x,y)\bar\pi(dxdy)
+\int_{\partial B_R\times B_R}f_3^2(z,y)\bar\pi(dzdy)\nonumber\\
&+\int_{B_R\times \partial B_R}f_4^2(x,w)\bar\pi(dxdw)
+\int_{\partial B_R\times \partial B_R}f_5^2(z,w)\bar\pi(dzdw).
\end{align}
%
By the triangle inequality w.~r.~t.~$L^2(\lambda_z)$ and $L^2(\mu_w)$,
and using that $\lambda_z$, $\mu_w$ are probability measures, we obtain
\begin{align}\label{ao82}
f_3\le|z-y|+\tilde f_3,\quad
f_4\le|x-w|+\tilde f_4,\quad
f_5\le|z-w|+\sqrt{2}\tilde f_5,
\end{align}
where the three functions $\tilde f_3,\tilde f_4,\tilde f_5\ge 0$ are defined by
\begin{align}
\tilde f_3^2(z,y)&:=\int|x-z|^2\lambda_z(dx),\quad
\tilde f_4^2(x,w):=\int|w-y|^2\mu_w(dy),\label{ao78}\\
\tilde f_5^2(z,w)&:=\tilde f_3^2(z,y)+\tilde f_4^2(x,w).\label{ao77}
\end{align}
The factor of $\sqrt{2}$ in (\ref{ao82}) arises because of $\tilde f_3+\tilde f_4
\le\sqrt{2}\tilde f_5$.

\medskip

From (\ref{ao82}) we obtain by the triangle inequality for $\|\cdot\|$
\begin{align}\label{ao76}
\lefteqn{\big(\int_\Omega|x-y|^2d\pi\big)^\frac{1}{2}
\stackrel{(\ref{ao81})}{\le}\|(f_2,f_3,f_4,f_5)\|}\nonumber\\
&\le\|(|x-y|,|z-y|,|x-w|,|z-w|)\|
+\sqrt{2}\|(0,\tilde f_3,\tilde f_4,\tilde f_5)\|,
\end{align}
where we gave up a factor of $\sqrt{2}$ on $\tilde f_3,\tilde f_4$.
By definition (\ref{ao75}) of $\|\cdot\|$, the first r.~h.~s.~term in (\ref{ao76})
coincides with the square root of $\int_{\bar B_R\times\bar B_R}|x$ $-y|^2\bar\pi$,
which by optimality of $\bar\pi$ is $W(\lambda\llcorner B_R+f,$ $\mu\llcorner B_R+g)$,
as desired.

\medskip

Using the definitions (\ref{ao75}), (\ref{ao78}), and (\ref{ao77}), 
the square of the second r.~h.~s.~term in (\ref{ao76}) is 
equal to $2$ times
$\int_{\partial B_R\times\mathbb{R}^d}$ $\int|x$ $-z|^2\lambda_z(dx)\bar\pi(dzdy)$
$+\int_{\mathbb{R}^d\times\partial B_R}$ $\int|w-y|^2\mu_y(dy)\bar\pi(dxdw)$.
By admissibility of $\bar\pi$, this sum is equal to
$\int_{\partial B_R}\int|x-z|^2\lambda_z(dx)f(dz)$
$+\int_{\partial B_R}\int|w-y|^2\mu_y(dy)g(dw)$.
By definitions (\ref{ao79}) and (\ref{ao80}), this coincides with
$\int_{\Omega\cap\{X(\sigma)\in\partial B_R\}}|x-X(\sigma)|^2d\pi$
$+\int_{\Omega\cap\{X(\tau)\in\partial B_R\}}|X(\tau)-y|^2d\pi$.
Since 
we have $|x-X(\sigma)|^2$ $+|X(\tau)-y|^2$ $\le|x-y|^2$, this sum is $\le$ 
$\int_{\Omega\cap(\{X(\sigma)\in\partial B_R\}\cup
\{X(\tau)\in\partial B_R\})}|x-y|^2d\pi$. Note that this set of integration
coincides with $\Omega\cap\{\exists t\in[0,1]\;X(t)\in\partial B_R\}$, as desired.


\subsection{Constructing a competitor based on the Neumann-Poisson problem}
As mentioned after Corollary \ref{cor:orth}, the remaining task is to estimate
the first r.~h.~s.~line of (\ref{ao10}). For this, we will use Lemma \ref{lem:opt}
and construct a competitor for $W(\lambda\llcorner B_R+f,\mu\llcorner B_R+g)$
based on $\phi$,
%
%
the solution of the Neumann-Poisson problem (\ref{ao12}),
where we momentarily
think of the measures $\lambda,\mu$ as having continuous densities with respect
to the Lebesgue measure.

\begin{lemma}\label{lem:comp2}
\begin{align}\label{ao20}
W^2(\lambda\llcorner B_R+f,\mu\llcorner B_R+g)
\le\frac{1}{\min\{\min_{\bar B_R}\lambda,\min_{\bar B_R}\mu\}}\int_{B_R}|\nabla\phi|^2.
\end{align}
\end{lemma}

Lemma \ref{lem:comp2} makes a second dilemma apparent: The intention was to use
it in conjunction with Lemma \ref{lem:opt} to obtain an estimate on the first
r.~h.~s.~line in (\ref{ao10}). This however would require that we have
${\min_{\bar B_R}\lambda,\min_{\bar B_R}\mu\gtrapprox 1}$, so a (one-sided) closeness 
of $\mu$ and $\lambda$ to the Lebesgue measure in a strong topology,
as opposed to the closeness in a weak topology as expressed by (\ref{ao88}). 
Hence this provides another reason for approximating $\lambda$ and $\mu$
by more regular versions.

\medskip

{\sc Proof Lemma \ref{lem:comp2}}. The proof is short if one uses the Benamou-Brenier
formulation in its distributional version, as we shall do. We recommend 
\cite[Section 6.1]{Santambrogio} to the reader regarding more details 
on the Benamou-Brenier formulation.
For every $t\in[0,1]$ we introduce the (singular non-negative) measure 
\begin{align}\label{ao84}
\rho_t:=t(\mu\llcorner B_R+g)+(1-t)(\lambda\llcorner B_R+f)
\end{align}
and the ($t$-independent) vector-valued measure
\begin{align}\label{ao85}
j_t:=\nabla\phi dx\llcorner B_R.
\end{align}
We note that (\ref{ao12}) in its distributional form of (\ref{ao62}) can be re-expressed as
\begin{align}\label{ao83}
\frac{d}{dt}\int\zeta d\rho_t=\int\nabla\zeta\cdot dj_t
\end{align}
for all test functions $\zeta$. In the jargon of the Benamou-Brenier formulation,
which is inspired from continuum mechanics, $\rho_t$ is a (mass) density, $j_t$ is a flux,
and (\ref{ao83}) is the distributional version of the continuity equation
$\partial_t\rho_t+\nabla\cdot j_t=0$
expressing conservation of mass.

\medskip

Following Benamou-Brenier one takes the Radon-Nikodym
derivative $\frac{dj_t}{d\rho_t}$ of the (vectorial) measure $j_t$ w.~r.~t.~$\rho_t$
(it plays the role of an Eulerian velocity field),
and considers the expression that corresponds
to the total kinetic energy:
\begin{align}\label{ao86}
\frac{1}{2}\int\bigg|\frac{dj_t}{d\rho_t}\bigg|^2 d\rho_t:
=\sup\bigg\{\int\xi\cdot dj_t-\int\frac{1}{2}|\xi|^2d\rho_t\bigg\}\in[0,\infty],
\end{align}
where the supremum is taken over all continuous vector fields $\xi$ with compact support.
Benamou-Brenier (see \cite[Section 5.4]{Santambrogio}) gives
\begin{align}\label{ao87}
W^2(\rho_0,\rho_1)\le\int_0^1\int\bigg|\frac{dj_t}{d\rho_t}\bigg|^2 d\rho_t dt.
\end{align}

\medskip

Since in our case, $j_t$ is supported in (the open) $B_R$, see (\ref{ao85}),
in the r.~h.~s.~of
(\ref{ao86}) we may restrict ourselves to $\xi$ supported in $B_R$. For these $\xi$'s,
definition (\ref{ao84}) yields $\int\xi\cdot dj_t-\int\frac{1}{2}|\xi|^2d\rho_t$ 
$=\int_{B_R}\big(\xi\cdot\nabla\phi$ $-\frac{1}{2}|\xi|^2(t\mu+(1-t)\lambda)\big)$. 
By Young's inequality in form of $\xi\cdot\nabla\phi$ 
$\le\frac{1}{2}(t\mu+(1-t)\lambda)|\xi|^2+\frac{1}{2(t\mu+(1-t)\lambda)}|\nabla\phi|^2$
we thus obtain for the r.~h.~s.~of (\ref{ao87})
\begin{align*}
\int\bigg|\frac{dj_t}{d\rho_t}\bigg|^2 d\rho_t\le\int_{B_R}\frac{|\nabla\phi|^2}{t\mu+(1-t)\lambda}
\le\frac{1}{\min\{\min_{\bar B_R}\lambda,\min_{\bar B_R}\mu\}}\int_{B_R}|\nabla\phi|^2.
\end{align*}
Since by definition (\ref{ao84}), the l.~h.~s.~of (\ref{ao87}) coincides
with the l.~h.~s.~of (\ref{ao20}), we are done.


\section{Harmonic approximation}\label{sec:harm}

The purpose of this section is to establish that the displacement in an optimal
plan $\pi$ can locally be approximated by a harmonic gradient $\nabla\phi$ (by which we
mean that for each Cartesian direction $i=1,\cdots,d$, 
the component $\partial_i\phi$ is harmonic, as 
a consequence of $-\triangle\phi=const$). This holds provided we are in the
perturbative regime, see Subsection \ref{ss:pert}, where $E$ and $D$ are defined.
More precisely, given any fraction $0<\theta\ll 1$, there exists a threshold
$\epsilon>0$ for $E+D$ so that below that threshold, the l.~h.~s.~of (\ref{ao89})
is only a fraction $\theta$ of $E$, plus a possibly large multiple of $D$.

\begin{proposition}\label{prop:harmonicApproximation}
For every $\theta>0$, there exist $\epsilon(d,\theta)>0$ and $C(d,\theta)<\infty$
such that the following holds. Let $\pi$ be optimal for $\lambda,\mu$;
provided $E+D\le\epsilon$, there exists a harmonic $\nabla\phi$ on $B_1$ such that
\begin{align}
\int_{(B_1\times\mathbb{R}^d)\cup(\mathbb{R}^d\times B_1)}
|(y-x)-\nabla\phi(x)|^2d\pi\le\theta E+CD,\label{ao89}\\
\int_{B_1}|\nabla\phi|^2\le C(E+D).\label{ao93}
\end{align}
\end{proposition}

(The proof actually reveals an explicit dependence of $\epsilon$ and $C$ on $\theta$.)
We will obtain $\nabla\phi$ by solving the Neumann-Poisson problem
\begin{align}\label{ao92}
-\triangle\phi=\frac{\mu(B_R)}{|B_R|}-\frac{\lambda(B_R)}{|B_R|}\;\mbox{in}\;B_R\quad\mbox{and}\quad
\nu\cdot\nabla\phi=\bar g-\bar f\;\mbox{on}\;\partial B_R,
\end{align}
where $\bar f,\bar g$ are suitable regular approximations of $f,g$, 
which are constructed in the nonlinear approximation Lemma \ref{lem:approx}, 
which also guides the choice of $R\in[2,3]$. (In fact, in Subsection \ref{PDEestimates}, we
will replace $\bar g-\bar f$ by its mollification.)
We note that $\bar f,\bar g\in L^2(\partial B_R)$ provides sufficient
regularity: Indeed, according to (\ref{ao48}) and (\ref{ao92}) we have $\int_{B_R}|\nabla\phi|^2$
$=\int_{\partial B_R}\phi (\bar g-\bar f)$, recalling the normalization $\int_{B_R}\phi=0$.
Applying Cauchy-Schwarz and then the Poincar\'e-trace estimate
$\int_{\partial B_R}\phi^2$ $\le C_P\int_{B_R}|\nabla\phi|^2$, we obtain
\begin{align}\label{as24}
\int_{B_R}|\nabla\phi|^2\le C_P\int_{\partial B_R}(\bar g-\bar f)^2.
\end{align}
In particular, \eqref{ao93} is a consequence of (\ref{ao96}),
for a suitable choice of $R\in[2,3]$.

\medskip

By an application of Lemma \ref{lem:comp2} to the setting of (\ref{ao92}), we have
\begin{align}\label{ao90}
W^2(\frac{\lambda(B_R)}{|B_R|}dx\llcorner B_R+\bar f,&
\frac{\mu(B_R)}{|B_R|}dx\llcorner B_R+\bar g)\nonumber\\
&\le\frac{1}{\min\{\frac{\lambda(B_R)}{|B_R|},\frac{\mu(B_R)}{|B_R|}\}}\int_{B_R}|\nabla\phi|^2.
\end{align}
Working with (\ref{ao92}) instead of (\ref{ao12})
creates the additional task of estimating the first r.~h.~s.~term (\ref{ao70}) 
of Lemma \ref{lem:opt} by the l.~h.~s.~of (\ref{ao90}), which is conveniently done
with help of the triangle inequality:
\begin{align}\label{ao91}
W(\lambda\llcorner B_R+f,&\mu\llcorner B_R+g)
\le W(\frac{\lambda(B_R)}{|B_R|}dx\llcorner B_R+\bar f,
\frac{\mu(B_R)}{|B_R|}dx\llcorner B_R+\bar g)\nonumber\\
&+W(\lambda\llcorner B_R,\frac{\lambda(B_R)}{|B_R|}dx\llcorner B_R)
+\mbox{same term with $\lambda\leadsto \mu$}\nonumber\\
&+W(f,\bar f)+\mbox{same term with $f\leadsto g$}.
\end{align}
We now return to the first r.~h.~s.~line in Corollary \ref{cor:orth}; in view of the
elementary
\begin{align}\label{cw05}
\lefteqn{\int_\Omega|x-y|^2d\pi-\int_{B_R}|\nabla\phi|^2}\nonumber\\
&\le 2\big(\int_\Omega|x-y|^2d\pi\big)^\frac{1}{2}
\Big(\big(\int_\Omega|x-y|^2d\pi\big)^\frac{1}{2}
-\big(\int_{B_R}|\nabla\phi|^2\big)^\frac{1}{2}\Big),
\end{align}
and noting that by definitions (\ref{ao17}) and (\ref{ao45}),
the first r.~h.~s.~factor is estimated by $E$; by Young's inequality,
it suffices to estimate the second factor.
Combining (\ref{ao70}), (\ref{ao90}) and (\ref{ao91}) we see that 
it is $\le$
\begin{align}\label{ao94}
\big(\frac{1}{(\min\{\frac{\lambda(B_R)}{|B_R|},\frac{\mu(B_R)}{|B_R|}\})^\frac{1}{2}}-1\big)
\big(\int_{B_R}|\nabla\phi|^2\big)^\frac{1}{2}
\end{align}
plus 
\begin{align}\label{ao95}
\lefteqn{\big(2\int_{\Omega\cap\{\exists t\in[0,1]\;X(t)\in\partial B_R\}}|x-y|^2d\pi\big)^\frac{1}{2}}
\nonumber\\
&+W(\lambda\llcorner B_R,\frac{\lambda(B_R)}{|B_R|}dx\llcorner B_R)
+\mbox{same term with $\lambda\leadsto \mu$}\nonumber\\
&+W(f,\bar f)+\mbox{same term with $f\leadsto g$}.
\end{align}
We expect (and will show for a good radius $R\in[2,3]$) 
that in the regime $D=o(1)$, the prefactor on the r.~h.~s.~of (\ref{ao94})
is $O(\sqrt D)$, and that the second line in (\ref{ao95})
is $O(\sqrt{D})$, as consistent with (\ref{ao89}). These two technicalities 
are stated in Lemma \ref{lem:localiseD}. The main task is thus to control
the last line in (\ref{ao95}).


\subsection{Approximating the boundary data}

The main remaining task is to identify a good radius $R\in[2,3]$ 
and to construct $\bar f$ and $\bar g$. Again, there is a trade-off/conflict of interest:
\begin{itemize}
\item On the one hand, the Neumann boundary data $\bar g-\bar f$ have to be sufficiently regular so
that the solution $\phi$ of (\ref{ao92}) is. In particular, we need 
(\ref{ao93}) (with $B_1$ replaced by the larger $B_R$) to obtain that the error 
(\ref{ao94}) is $o(E+D)$. Via (\ref{as24}), this is ensured
by (\ref{ao96}) in the upcoming Lemma \ref{lem:approx}. In fact, it even yields
uniform integrability of $|\nabla\phi|^2$ on $B_R$, which is crucial to show that
also the last line in (\ref{ao10}) is $o(E+D)$. 
\item On the other hand, $(\bar f,\bar g)$ has to be sufficiently close to
$(f,g)$. In particular, in view of the last term in (\ref{ao95}) we need 
$W^2(f,\bar f)$ $+W^2(g,\bar g)$ $=o(E)+O(D)$. This is ensured
by (\ref{ao97}) 
in the upcoming Lemma \ref{lem:approx}. Here, as for (\ref{ao70}), we 
will eventually need to appeal to
$\int_1^2\int_{\{\exists t\in[0,1]\;X(t)\in\partial B_R\}}|x-y|^2d\pi\,dR$ $=o(E+D)$,
see Subsection \ref{sec:crossing}.
\end{itemize}
In the upcoming approximation lemma we restrict to $g$ for brevity.


\begin{lemma}\label{lem:approx}
We suppose that 
\begin{align}\label{as22}
(X\in\Omega \;\Longrightarrow\; y\in B_5)
\quad\mbox{for}\;(x,y)\in{\rm supp}\pi.
\end{align}
Then for every $R\in[2,3]$ there exists a non-negative function $\bar g_R$ on $\partial B_R$ 
such that 
\begin{align}
W^2(g_R,\bar g_R)&\le 8\big(
\int_{\Omega\cap\{\exists t\in[0,1]\;X(t)\in\partial B_R\}}|x-y|^2d\pi
+D\big),\label{ao97}\\
\int_2^3 \int_{\partial B_R}\bar g_R^2\,dR&\le 5^{d-1}\kappa_\mu(3E+D).\label{ao96}
\end{align}
\end{lemma}

Note that we put an index $R$ on $g$ because the definition (\ref{ao09bis}) 
obviously depends on $R$. 

\medskip

{\sc Proof of Lemma \ref{lem:approx}}.
We fix an $R\in[2,3]$ and start with the construction of $\bar g_R$, momentarily
returning to our short-hand notation $\bar g$.
Let $\bar\pi$ be optimal for $W^2(\mu\llcorner B_5,\kappa_\mu dz\llcorner B_5)$;
note that $\bar\pi$ is supported on $B_5\times B_5$. We extend it (trivially)
by the identity to $\mathbb{R}^d\times\mathbb{R}^d$; the extension (which we still call) 
$\bar\pi$ is admissible for $W^2(\mu,\kappa_\mu dz\llcorner B_5+\mu\llcorner B_5^c)$. We retain
\begin{align}\label{ho05}
\int|y-z|^2d\bar\pi= W^2(\mu\llcorner B_5,\kappa_\mu dz\llcorner B_5)
\stackrel{(\ref{ao88})}{\le} D.
\end{align}
Like in the proof of the triangle inequality for the Wasserstein metric, 
we disintegrate $\bar\pi$ according to
\begin{align}\label{ho09}
\int\zeta(y,z)\bar\pi(dz|y)\mu(dy)=\int\zeta(y,z)\bar\pi(dydz),
\end{align}
since this family of (conditional) probability measures 
$\{\bar\pi(\cdot|y)\}_{y\in\mathbb{R}^d}$ allows us to define
the measure $\tilde\pi$ on $\mathbb{R}^d\times\mathbb{R}^d\times\mathbb{R}^d$ via
\begin{align*}
\int\zeta(x,y,z)\tilde\pi(dxdydz)=\int\int\zeta(x,y,z)\bar\pi(dz|y)\pi(dxdy),
\end{align*}
which has the desired property
\begin{align}\label{ao99}
\begin{array}{cc}
\mbox{marginal of}\;\tilde\pi\;\mbox{w.~r.~t.}\;(x,y)&\;\;=\;\;\pi,\\
\mbox{marginal of}\;\tilde\pi\;\mbox{w.~r.~t.}\;(y,z)&\;\;=\;\;\bar\pi.
\end{array}
\end{align}

\medskip

\begin{figure}[h!]
  \centering
    \includegraphics[width=0.5\linewidth]{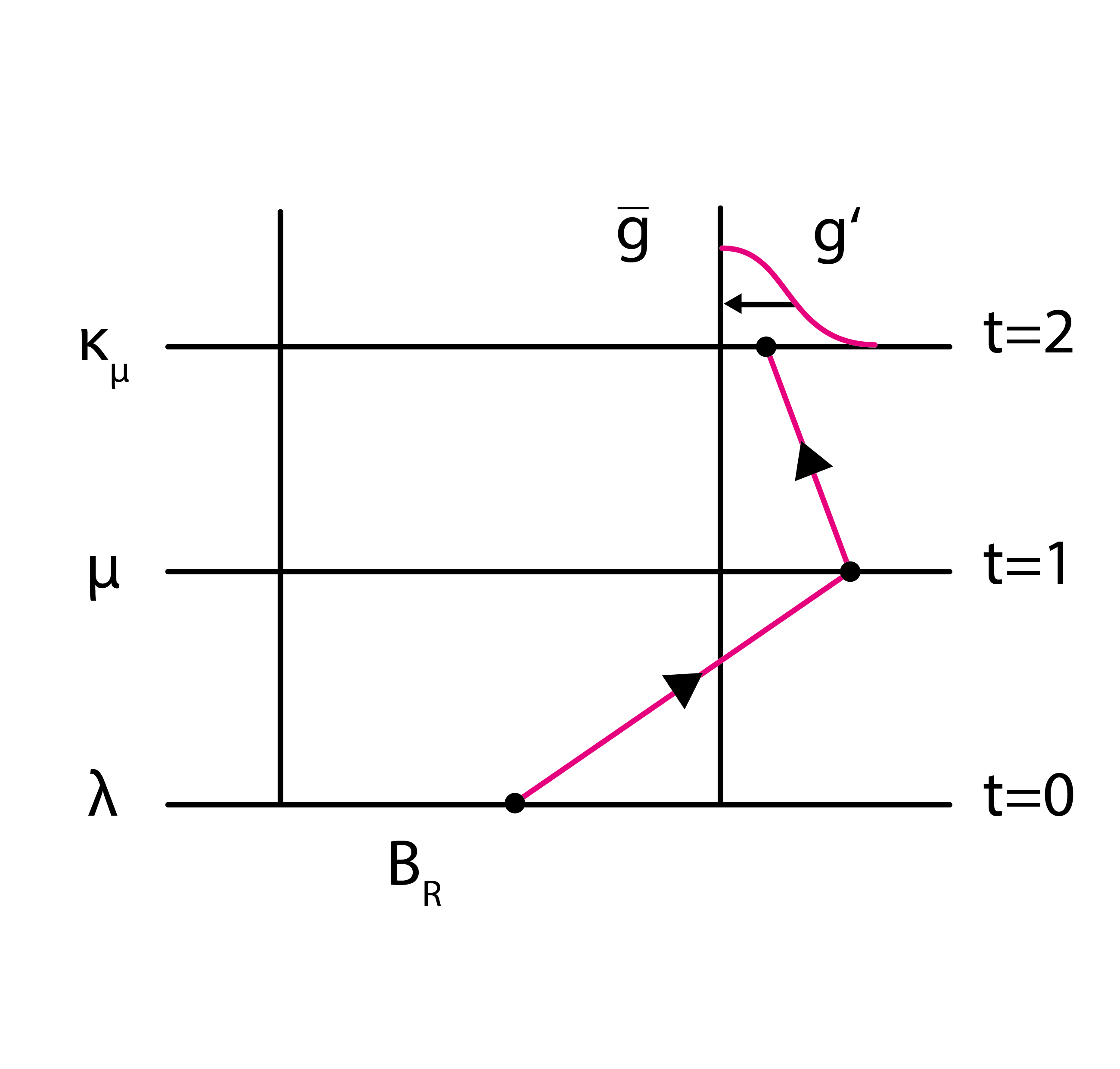}
    \caption{Defining $\bar g$ and $g'$}
\label{fig:gs}
\end{figure}

We still associate the trajectory $X$ to a pair $(x,y)$ as in (\ref{ao00}).
We extend this definition by associating to a triplet $(x,y,z)$ the continuous
piecewise affine trajectory defined by (\ref{ao00}) followed by
\begin{align*}
X(t)=(t-1)z+(2-t)y\quad\mbox{for}\;t\in[1,2].
\end{align*}
We now are interested in the distribution $g'$ of the endpoint $z=X(2)$ of those
trajectories that exit $\bar B_R$ during the first time interval $[0,1]$,
i.~e.~those trajectories $X\in\Omega$ with $X(\tau)\in\partial B_R$, see Figure \ref{fig:gs}. This distribution
is defined by
\begin{align}\label{ho01}
\int\zeta d g'=\int_{\Omega\cap\{X(\tau)\in\partial B_R\}}\zeta(z)\tilde\pi(dxdydz).
\end{align}

We learn from (\ref{as22}) that
$y=X(1)\in B_5$ for the trajectories $X$ that contribute to (\ref{ho01}).
This yields $z=X(2)\in B_5$ by
$\bar\pi(B_5,B_5^c)=0$ according to our extension of $\bar\pi$. Hence we have
\begin{align}\label{ho14}
g'\;\mbox{is supported in}\;B_5.
\end{align}
Therefore by the second item
in (\ref{ao99}) and the admissibility of $\bar\pi$ for 
$W^2(\mu,$ $\kappa_\mu dz\llcorner B_5$ $+\mu\llcorner B_5^c)$, we obtain
\begin{align*}
\int\zeta d g'
\le\int_{\{z\in B_5\}}\zeta(z)\bar\pi(dydz)=\kappa_\mu\int_{B_5}\zeta
\quad\mbox{provided}\;\zeta\ge 0.
\end{align*}
Hence $g'$ admits a Lebesgue density, which we still denote by $g'$ and that satisfies
\begin{align}\label{ho11}
g'\le\kappa_\mu.
\end{align}
%
%
Finally, we radially project $g'$ onto $\partial B_R$:
\begin{align}\label{ho02}
\int\zeta d\bar g=\int\zeta(R\frac{z}{|z|})g'(dz),
\end{align}
see Fig. \ref{fig:gs}.
This concludes the construction of $\bar g$, we now turn to its estimate.

\medskip

We start with (\ref{ao97}) and note that an admissible plan for $W^2(g,\bar g)$ is given by
\begin{align*}
\int_{\Omega\cap\{X(\tau)\in\partial B_R\}}\zeta(X(\tau),R\frac{z}{|z|})d\tilde\pi.
\end{align*}
Indeed, on the one hand,
for $\zeta$ only depending on the first variable, $\tilde\pi$ may be replaced by $\pi$
according to the first item in (\ref{ao99})
so that we obtain $\int\zeta dg$ by its definition (\ref{ao09bis}).
On the other hand, for $\zeta$ only depending on the second variable,
we obtain $\int\zeta d\bar g$ by combining (\ref{ho01}) and (\ref{ho02}).
Hence we have
\begin{align}\label{ho06}
W^2(g,\bar g)\le\int_{\Omega\cap\{X(\tau)\in\partial B_R\}}|X(\tau)-R\frac{z}{|z|}|^2d\tilde\pi.
\end{align}
Since $X(\tau)\in\partial B_R$, an elementary geometric argument on the radial projection yields
for the integrand
$|X(\tau)-R\frac{z}{|z|}|$ $\le2|X(\tau)-z|$, so that
$|X(\tau)-R\frac{z}{|z|}|^2$ $\le8(|x-y|^2+|y-z|^2)$.
This allows us to appeal to the compatibility (\ref{ao99}) 
and to (\ref{ho05}):
\begin{align*}
W^2(g,\bar g)\le 8\big(\int_{\Omega\cap\{\exists t\in[0,1]\;X(t)\in\partial B_R\}}
|x-y|^2d\pi+D\big).
\end{align*}
\ignore{BEGIN IGNORE
We also note that for $X$ in the domain of integration $\Omega\cap\{X(\tau)\in\partial B_R\}$ 
we have $\min_{[0,1]}|X|\le R\le\max_{[0,1]}|X|$; by $R\le 2$ and definition (\ref{ho03}),
this implies $\min_{[0,1]}|X|\le R\le\min_{[0,1]}|X|+M$ on the support of $\tilde\pi$. 
We recall that by $M<1$ and $R\le 2$ we also have $y\in B_3$ there.
Summing up, we infer from (\ref{ho06})
\begin{align*}
\lefteqn{W^2(g,\bar g)}\nonumber\\
&\le4\int_{\{y\in B_3\}\cap
\{\min_{[0,1]}|X|\le R\le\min_{[0,1]}|X|+M\}}(|x-y|+|y-z|)^2d\tilde\pi.
\end{align*}
We now (re)introduce the index $R$ and integrate over $R\in[1,2]$;
expressing the domain of integration $\{\min_{[0,1]}|X|\le R\le\min_{[0,1]}|X|+M\}$
by a characteristic function and 
exchanging the order of integration we obtain
\begin{align*}
\int_1^2 W^2(g_R,\bar g_R)dR\le 4M\int_{\{y\in B_3\}}(|x-y|+|y-z|)^2d\tilde\pi.
\end{align*}
The use of $(|x-y|+|y-z|)^2$ $\le 2(|x-y|^2+|y-z|^2)$
allows us to appeal to the compatibility (\ref{ao99}):
\begin{align*}
\int_1^2 W^2(g_R,\bar g_R)\le 8 M\big(\int_{\mathbb{R}^d\times B_3}
|x-y|^2d\pi+\int|y-z|^2d\bar\pi\big).
\end{align*}
By definition (\ref{ao45}) and by (\ref{ho05}) this turns into (\ref{ao97}).
}

\medskip

In preparation for establishing (\ref{ao96}), we first provide an estimate
of the measure $g'$ defined in (\ref{ho01}), which shows that it is concentrated
near $\partial B_R$, see (\ref{ho07}). By definition (\ref{ho01}) we have
\begin{align*}
\int||z|-R|dg'=\int_{\Omega\cap\{X(\tau)\in\partial B_R\}}||z|-R|\tilde\pi(dxdydz).
\end{align*}
Since $|X(\tau)|=R$, we may write $||z|-R|$ $=||z|-|X(\tau)||$ $\le|x-y|$ $+|y-z|$.
Since $X(\tau)\in\partial B_R$ implies $\min_{[0,1]}|X|\le R\le\max_{[0,1]}|X|$,
and $X(1)\in B_5$ by assumption (\ref{as22}), we thus obtain
\begin{align*}
\lefteqn{\int||z|-R|dg'}\nonumber\\
&\le\int_{\{y\in B_5\}\cap\{\min_{[0,1]}|X|\le R\le\max_{[0,1]}|X|\}}
(|x-y|+|y-z|)\tilde\pi(dxdydz).
\end{align*}
Making the index $R$ appear and integrating in $R$, this gives
\begin{align*}
\lefteqn{\int_2^3\int||z|-R|dg'_R\,dR}\nonumber\\
&\leq \int_{\{y\in B_5\}}(\max_{[0,1]}|X|-\min_{[0,1]}|X|)
(|x-y|+|y-z|)\tilde\pi(dxdydz).
\end{align*}
Using $\max_{[0,1]}|X|-\min_{[0,1]}|X|$ $\le|x-y|$ and then Young's inequality
in form of $|x-y|(|x-y|+|y-z|)$ $\le\frac{3}{2}|x-y|^2+\frac{1}{2}|y-z|^2$ we 
thus obtain from (\ref{ao99})
\begin{align*}
\lefteqn{\int_2^3\int||z|-R|dg'_R\,dR}\nonumber\\
&=\frac{3}{2}\int_{\{y\in B_5\}}|x-y|^2\pi(dxdy)
+\frac{1}{2}\int|y-z|^2\bar\pi(dydz).
\end{align*}
By definition (\ref{ao45}) and by (\ref{ho05}) this turns into 
\begin{align}\label{ho07}
\int_2^3\int||z|-R|dg'_R\,dR\le\frac{3}{2}E+\frac{1}{2}D.
\end{align}

\medskip

In order to pass from (\ref{ho07}) to (\ref{ao96}) we need
\begin{align}\label{ho08}
\int_{\partial B_R}\frac{1}{2}\bar g^2\le 5^{d-1}\kappa_\mu\int||z|-R|dg'.
\end{align}
Here comes the argument for (\ref{ho08}):
By (\ref{ho11}), it follows once we show for some density $g'$ with (\ref{ho14})
\begin{align}\label{ho12}
\int_{\partial B_R}\frac{1}{2}\bar g^2\le 5^{d-1}({\rm ess\, sup\,}g') 
\int||z|-R|g' dz.
\end{align}
Introducing polar coordinates $z=r\hat z$ with $r\in(0,\infty)$ and $\hat z\in\partial B_1$, 
which are natural to re-express (\ref{ho02}), (\ref{ho12}) reduces to
the single-variable statement
\begin{align*}
\frac{1}{2}\big(\int g' r^{d-1}dr\big)^2
\le 5^{d-1}({\rm ess \, sup\,}\kappa) \int|r-R|g' r^{d-1}dr.
\end{align*}
It is convenient to rephrase this in terms of $\tilde g=g' r^{d-1}$;
since because of (\ref{ho14}) we have ${\rm ess \, sup\,}\tilde g$
$\le 5^{d-1}{\rm ess \, sup\,}g'$, it suffices to show for an arbitrary
function $\tilde g\ge 0$ of $r\in(-\infty,\infty)$ that
\begin{align}\label{ho13}
\frac{1}{2}\big(\int\tilde g dr\big)^2
\le ({\rm ess \, sup\,}\tilde g) \int|r-R|\tilde g dr.
\end{align}
The argument for (\ref{ho13}) is elementary: By translation in $r$,
we may assume $R=0$; by homogeneity in $\tilde g$, we may assume
${\rm ess \, sup\,}\tilde g=1$, that is, $\tilde g\in[0,1]$.
We now change perspective and seek to minimize the r.~h.~s.~$\int|r|\tilde g dr$ under
constraining the l.~h.~s.~ through prescribing $m=\int\tilde g dr$.
By linearity of $\int|r|\tilde g dr$ in $\tilde g$, this functional
assumes its minimum on extremal points w.~r.~t.~the constraints
$\tilde g\in[0,1]$ and $\int\tilde g dr=m$. 
Those are characteristic functions of sets
of Lebesgue measure $m$. Clearly, the set $I$ with $|I|=m$ that minimizes 
$\int_I|r|dr$ is given by $I=[-\frac{m}{2},\frac{m}{2}]$;
the minimum is $\frac{m^2}{2}$, as desired.


\subsection{Crossing trajectories}\label{sec:crossing}

In view of the second r.~h.~s.~term in (\ref{ao70}) of Lemma \ref{lem:opt} and
the first r.~h.~s.~term in (\ref{ao97}) of Lemma \ref{lem:approx},
we need to argue that for a suitable radius $R$, the trajectories crossing $\partial B_R$
do not contribute much to $E$. This is the only part of the argument where
we directly rely on the optimality criterion for $\pi$, namely the cyclical
monotonicity of its support, c.~f. \cite[Section 1.6.2]{Santambrogio}.
We just need it in form of plain monotonicity:
\begin{align}\label{m08}
(x-x')\cdot(y-y')\geq 0\quad\mbox{for all}\;(x,y), (x',y')\in{\rm supp}\pi.
\end{align}
This in fact implies that all trajectories are short
in our regime of $E+D\ll 1$:

\begin{lemma}\label{lem:linfty}
We have
\begin{align}
|x-y|=o(1)\quad\mbox{for}\;(x,y)
\in\big((B_4\times\mathbb{R}^d)\cup(\mathbb{R}^d\times B_4)\big)\cap{\rm supp}\pi,\label{m07}\\
\int_2^3\int_{\Omega\cap\{\exists t\in[0,1]\; X(t)\in \partial B_R\}}|x-y|^2 d\pi dR
=o(E+D),\label{m05}\\
\int_2^3\pi\big(\Omega\cap\{\exists t\in[0,1]\; X(t)\in \partial B_R\}\big)dR=o(1).\label{m06}
\end{align}
\end{lemma}

As a consequence of (\ref{m05}) we may indeed chose $R\in[2,3]$ 
so that the terms in \eqref{ao70} and \eqref{ao97} are $o(E)+O(D)$.
As a consequence of (\ref{m07}), also \eqref{as22} is satisfied.

\medskip

{\sc Proof of Lemma \ref{lem:linfty}}. 
We start by deriving (\ref{m05}) and (\ref{m06}) from (\ref{m07}).
As in the proof of (\ref{ao96}) we note that $X(t)\in\partial B_R$
implies $\min_{[0,1]}|X|\le R\le \max_{[0,1]}|X|$ and thus
\begin{align*}
\lefteqn{\int_2^3\int_{\Omega\cap\{\exists t\in
[0,1]\;X(t)\in\partial B_R\}}
|x-y|^2\pi(dxdy)dR}\nonumber\\
&\le\int_{(B_3\times \R^d)\cup (\R^d\times B_3)}(\max_{[0,1]}|X|-\min_{[0,1]}|X|)|x-y|^2\pi(dxdy).
\end{align*}
Since $\max_{[0,1]}|X|-\min_{[0,1]}|X|$ $\le|x-y|$, (\ref{m05})
now follows from (\ref{m07}) by the definition (\ref{ao45}).
Likewise, we have 
\begin{align*}
&\int_2^3\pi(\Omega\cap\{\exists t\in[0,1]\; X(t)\in \partial B_R\}) dR\\
&=o(1)\pi((B_3\times \R^d)\cup (\R^d\times B_3)).
\end{align*}
By (\ref{ao00bis}) we obtain that $\pi((B_3\times \R^d)\cup (\R^d \times B_3))$ $\le\lambda(B_5)+\mu(B_5)$, so that we may
appeal to (\ref{ao88}) to the effect of $\lambda(B_5)$ $\le |B_5|(1+\sqrt{D})$ $\lesssim 1$. Here and in the following we use $\lesssim$ to mean up to constants that only depend on $d$.

\medskip

We now turn to proving \eqref{m07};
a more explicit but pedestrian argument can be found in \cite[Lemma 2.9]{Goldman2021}.
Since $E\ll 1$, we expect there to be many short trajectories. 
We exploit monotonicity in order to upgrade this to a statement about all trajectories. 
By definition (\ref{ao17}) and due to ($x\leftrightarrow y$)-symmetry, 
we may assume that $x\in B_4$.
We will use (\ref{m08}) in form of
\begin{align}\label{m02}
(x-y)\cdot(x-x')\leq \frac 3 2 |x-x'|^2+\frac 1 2 |x'-y'|^2.
\end{align}
For $\zeta\ge 0$ supported in $B_5$ we integrate 
\eqref{m02} against $\zeta(x')\pi(dx'dy')$.
Using the admissibility of $\pi$ and the definition (\ref{ao45}) of $E$, we find
\begin{align}\label{m04}
\lefteqn{(x-y)\cdot\int(x-x')\zeta(x')\lambda}(dx')\nonumber\\
&\le\frac{3}{2}\int |x-x'|^2\zeta(x')\lambda(dx')
+\frac{1}{2}(\sup\zeta) E.
\end{align}

\medskip

By definition (\ref{ao88}), $D\ll 1$ implies that $\lambda\llcorner B_5$ weakly
close to $dx'\llcorner B_5$. In order to obtain a uniformity when applying this
to (\ref{m04}),
we now specify a class of functions $\zeta$ in (\ref{m04}): Fixing a 
smooth mask $\hat\zeta\ge 0$ supported in $B_1$
and some radius $0<r\ll 1$ to be optimized later, 
we consider $\zeta(x')=\hat\zeta(\frac{x'-x_0'}{r})$ for any center $x_0'$,
only restricted by the requirement that $\zeta$ is supported in $B_5$.
This class of functions is compact with respect to the uniform topology so that we have
\begin{align*}
\int(x-x')\zeta(x')\lambda(dx')&=\int(x-x')\zeta(x')dx'+o_r(1),\\
\int |x-x'|^2\zeta(x')\lambda(dx')&=\int |x-x'|^2\zeta(x')dx'+o_r(1)
\end{align*}
as $D\ll 1$,
uniformly in the center but at fixed radius (as the subscript in $o_r(1)$ is to indicate). 
Specifying $\hat\zeta$ to have unit integral and vanishing first moment, this yields
\begin{align*}
\int(x-x')\zeta(x')\lambda(dx')&=r^d(x-x'_0)+o_r(1),\\
\int |x-x'|^2\zeta(x')\lambda(dx')&\lesssim r^d(|x-x_0'|^2+r^2)+o_r(1).
\end{align*}
This prompts the choice of the center $x_0'=x-r\frac{x-y}{|x-y|}$, 
which is admissible since in view of $x\in B_4$ we still have that $\zeta$ is supported
in $B_5$, to the effect of
\begin{align*}
(x-y)\cdot\int(x-x')\zeta(x')\lambda(dx')&=r^{d+1}|x-y|+o_r(|x-y|),\\
\int |x-x'|^2\zeta(x')\lambda(dx')&\lesssim r^{d+2}+o_r(1).
\end{align*}
Inserting this into (\ref{m04}) yields
\begin{align*}
|x-y|\lesssim r+\frac{E}{r^{d+1}}+o_r(|x-y|)+o_r(1).
\end{align*}

We now may conclude: We first choose $r$ so that the first
r.~h.~s.~term is small; we then choose $E$ so that the second term is small,
and we finally choose $D$ so small such that the third term may be absorbed
into the l.~h.~s.~ and such that the last term is small.


\subsection{Restricting the data term $D$}

While the data term $D$ is defined w.~r.~t.~to $B_5$, 
we rather need it w.~r.~t.~$B_R$ for our suitably chosen $R\in[2,3]$.
This is most prominent in the middle line of (\ref{ao95}) and the prefactor of (\ref{ao94}),
but also related to the middle and last line on the r.~h.~s.~of (\ref{ao10}),
as was discussed in Subsection \ref{ss:pert}.
Annoyingly, this restriction property for $D$ does not come for free and 
and requires arguments similar to the ones used in the proof of Lemma \ref{lem:approx}.
By symmetry, it is enough to consider $\lambda$.

\begin{lemma}\label{lem:localiseD} 
With $\kappa_R = \frac{\lambda(B_R)}{|B_R|}$ we have
\begin{align}\label{as21}
\int_2^3 W^2(\lambda\llcorner B_R,\kappa_R dx\llcorner B_R)+(\kappa_R-1)^2 dR=O(D).
\end{align}
\end{lemma}

For the proof of Lemma \ref{lem:localiseD}, it is convenient to have the following two
extensions of Lemma \ref{lem:comp2} on the relationship between OT and
the Poisson-Neumann problem; the first is a rough generalization,
the second provides the reverse relationship in a restricted setting:

\begin{corollary}\label{cor:2}
\begin{align}\label{as02}
W^2(\lambda\llcorner B_R+f,\mu\llcorner B_R+g)
\le\frac{4}{\min_{\bar B_R}\mu}\int_{B_R}|\nabla\phi|^2.
\end{align}
\end{corollary}

\begin{lemma}\label{lem:up}
Provided $f,g\equiv 0$ in (\ref{ao12}),
\begin{align}\label{as07}
W^2(\lambda\llcorner B_R,\mu\llcorner B_R)
\ge\frac{1}{\max\{\max_{\bar B_R}\lambda,\max_{\bar B_R}\mu\}}\int_{B_R}|\nabla\phi|^2.
\end{align}
\end{lemma}

{\sc Proof of Corollary \ref{cor:2}}. We first argue that for arbitrary $\lambda,\mu$
and $0\le M<\infty$,
\begin{align}\label{as03}
W(\lambda,\mu)\le\frac{1}{\sqrt{1+M}-\sqrt{M}}W(\lambda+M\mu,(1+M)\mu).
\end{align}
Indeed, by scaling we have $W(\lambda,\mu)$ $=\frac{1}{\sqrt{1+M}}W((1+M)\lambda,(1+M)\mu)$.
By the triangle inequality, $W((1+M)\lambda,(1+M)\mu)$
$\le W((1+M)\lambda,$ $\lambda+M\mu)$ $+W(\lambda+M\mu,(1+M)\mu)$, which we combine with
the obvious
$W((1+M)\lambda,\lambda+M\mu)$ $= W(M\lambda,M\mu)$.
Once more by scaling, $W(M\lambda,M\mu)$ $=\sqrt{M}W(\lambda,\mu)$, so that we may absorb.

\medskip

We now argue that with help of Lemma \ref{lem:comp2} we obtain 
\begin{align*}
W^2(\lambda\llcorner B_R+f,\mu\llcorner B_R+g)
\le\frac{1}{M(\sqrt{1+M}-\sqrt{M})^2}\frac{1}{\min_{\bar B_R}\mu}
\int_{B_R}|\nabla\phi|^2,
\end{align*}
which yields (\ref{as02}) under $M\uparrow\infty$. Indeed, we first apply
(\ref{as03}) with $(\lambda,\mu)$ replaced by $(\lambda\llcorner B_R+f,\mu\llcorner B_R+g)$.
We then use (\ref{ao20}) to estimate  
$W^2(\lambda\llcorner B_R+f+M(\mu\llcorner B_R+g),$ $(1+M)(\mu\llcorner B_R+g))$,
noting that when taking the difference of $(1+M)\mu$ and $\lambda+M\mu$
and of $(1+M)g$ and $f+Mg$, the $M$-dependent terms drop out and thus does not
affect the definition (\ref{ao12}) of $\phi$. It remains to observe that
the minimum of the Lebesgue densities of $(1+M)\mu$ and $\lambda+M\mu$ on $\bar B_R$
is bounded below by $M\min_{\bar B_R}\mu$.

\medskip

{\sc Proof of Lemma \ref{lem:up}}. We recall the Benamou-Brenier formulation
from Lemma \ref{lem:comp2}
and note that for the optimal $(\rho_t,j_t)$
\begin{align}\label{as05}
W^2(\lambda\llcorner B_R,\mu\llcorner B_R)\ge\int_0^1\int
\left|\frac{dj_t}{d\rho_t}\right|^2d\rho_tdt.
\end{align} 
Setting $M:=\max\{\max_{\bar B_R}\lambda,\max_{\bar B_R}\mu\}<\infty$, 
so that $\lambda\llcorner B_R,$ $\mu\llcorner B_R$ $\le Mdx\llcorner B_R$,
we have by McCann's displacement convexity (in conjunction with the convexity of $B_R$), c.f. \cite[Section 7.3]{Santambrogio},
\begin{align}\label{as04}
\rho_t\le Mdx\llcorner B_R
\end{align}
for all $t\in[0,1]$.
In view of definition (\ref{ao86}), we thus obtain for any test vector field $\xi$
\begin{align*}
\int\xi\cdot dj_t\le\frac{1}{2}\int|\frac{dj_t}{d\rho_t}|^2d\rho_t+\frac{M}{2}\int_{B_R}|\xi|^2.
\end{align*}
After integration in $t\in[0,1]$, replacing $\xi$ by a $a\xi$ and optimizing in
the constant $a>0$, this yields
$\int\xi\cdot dj_t$ $\le(\int|\frac{dj_t}{d\rho_t}|^2d\rho_t$ 
$M\int_{B_R}|\xi|^2)^\frac{1}{2}$.
In conjunction with (\ref{as05}) this implies
\begin{align}\label{as06}
\lefteqn{\int_0^1j_t dt\ll dx\llcorner B_R}\nonumber\\
&\mbox{and}\quad
\int_{B_R}|\int_0^1 j_t dt|^2\le MW^2(\lambda\llcorner B_R,\mu\llcorner B_R),
\end{align}
where we identify the measure $\int_0^1 j_t dt$ with its Lebesgue density. 

\medskip

From the admissibility condition (\ref{ao83}) we obtain by integration
in $t\in[0,1]$ and using $\rho_{t=0}=\lambda\llcorner B_R$, $\rho_{t=1}=\mu\llcorner B_R$
\begin{align*}
-\nabla\cdot\int_0^1j_tdt=\mu\llcorner B_R-\lambda\llcorner B_R\quad\mbox{distributionally
on}\;\mathbb{R}^d.
\end{align*}
%
Thus by (\ref{ao12}) (with $f,g\equiv0$), we have that
$\int_0^1j_tdt-\nabla\phi$ (with $\nabla\phi$ trivially extended beyond $\bar B_R$)
is distributionally divergence-free in $\mathbb{R}^d$. Extending $\phi$ in a
(compactly supported) $C^1$-manner outside of $\bar B_R$ and testing with this
extension, because $\int_0^1j_t dt$ is supported in $B_R$, cf.~(\ref{as06}), 
we obtain $\int_{B_R}\nabla\phi\cdot(\int_0^1j_tdt-\nabla\phi)$ $=0$. By Cauchy-Schwarz
this yields
\begin{align*}
\int_{B_R}|\nabla\phi|^2\le\int_{B_R}|\int_0^1j_tdt|^2.
\end{align*}
Combining this with (\ref{as06}) we obtain (\ref{as07}).

\medskip


{\sc Proof of Lemma \ref{lem:localiseD}}.
Inside this proof, we denote by $\pi$ the optimal plan in $W^2(\lambda\llcorner B_5,
\kappa_\lambda dx\llcorner B_5)$.
Following the proof of Lemma \ref{lem:approx}, we monitor where entering and exiting
trajectories end up:
\begin{align}
\int\zeta df'=\int_{\Omega\cap\{X(0)\not\in B_R\}\cap\{X(1)\in B_R\}}
\zeta(X(1))d\pi,\nonumber\\
\int\zeta dg'=\int_{\Omega\cap\{X(0)\in B_R\}\cap\{X(1)\not\in B_R\}}
\zeta(X(1))d\pi.\label{as18}
\end{align}

Clearly, the non-negative measures $f',g'$ are supported in $B_R$ and $\bar B_5-B_R$,
respectively, and satisfy $f',g'\le\kappa_\lambda$. We introduce the corresponding mass
densities w.~r.~t.~$B_R$
\begin{align}\label{as17}
\kappa_f:=\frac{f'(\mathbb{R}^d)}{|B_R|}\le\kappa_\lambda,
\quad\kappa_g:=\frac{g'(\mathbb{R}^d)}{|B_R|}.
\end{align}

\medskip

We start with the triangle inequality
\begin{align}\label{as10}
\lefteqn{W(\lambda\llcorner B_R,(\kappa_\lambda-\kappa_f+\kappa_g)dx\llcorner B_R)}\nonumber\\
&\le W(\lambda\llcorner B_R,\kappa_\lambda dx\llcorner B_R-f'+g')
\nonumber\\
&+W(\kappa_\lambda dx\llcorner B_R-f'+g',(\kappa_\lambda-\kappa_f)dx\llcorner B_R+g')
\nonumber\\
&+W((\kappa_\lambda-\kappa_f)dx\llcorner B_R+g',
(\kappa_\lambda-\kappa_f+\kappa_g)dx\llcorner B_R).
\end{align}
Restricting the optimal $\pi$ to trajectories that start in $B_R$,
we obtain an admissible plan for the first r.~h.~s.~term of (\ref{as10}):
\begin{align*}
W^2(\lambda\llcorner B_R,\kappa_\lambda dx\llcorner B_R-f'+g')
\le W^2(\lambda\llcorner B_5,\kappa_\lambda dx\llcorner B_5)\le D.
\end{align*}
We now turn to the second r.~h.~s.~term of (\ref{as10}),
and let $\phi'$ denote the solution of
\begin{align*}
-\triangle\phi'=f'-\kappa_f\;\;\mbox{in}\;B_R,\quad
\nu\cdot\nabla\phi'=0\;\;\mbox{on}\;\partial B_R.
\end{align*} 
On the one hand, we obtain from Corollary \ref{cor:2}
\begin{align*}
W(\kappa_\lambda dx\llcorner B_R-f'+g',(\kappa_\lambda-\kappa_f)dx\llcorner B_R+g')
\le\frac{4}{\kappa_\lambda-\kappa_f}\int_{B_R}|\nabla\phi'|^2.
\end{align*}
On the other hand, we obtain from Lemma \ref{lem:up} and $f'\le\kappa_\lambda$
\begin{align*}
\frac{1}{2\kappa_\lambda}\int_{B_R}|\nabla\phi'|^2
\le W(\kappa_\lambda dx\llcorner B_R,(\kappa_\lambda-\kappa_f)dx\llcorner B_R+f').
\end{align*}
The combination of these two inequalities yields
\begin{align}\label{as11}
\lefteqn{W(\kappa_\lambda dx\llcorner B_R-f'+g',
(\kappa_\lambda-\kappa_f)dx\llcorner B_R+g')}\nonumber\\
&\le \frac{8\kappa_\lambda}{\kappa_\lambda-\kappa_f}
W(\kappa_\lambda dx\llcorner B_R,(\kappa_\lambda-\kappa_f)dx\llcorner B_R+f').
\end{align}

\medskip

In view of (\ref{as11}), we may treat the second r.~h.~s.~term in
(\ref{as10}) alongside of the last term,
and focus on the latter. Like in Lemma \ref{lem:approx}, we introduce the projection
$\bar g$ of $g'$ onto $\partial B_R$, that is
\begin{align}\label{as15}
\int\zeta d\bar g=\int\zeta(\frac{Rx}{|x|})g'(dx).
\end{align}
We start with the triangle inequality
\begin{align}\label{as14}
\lefteqn{W((\kappa_\lambda-\kappa_f)dx\llcorner B_R+g',
(\kappa_\lambda-\kappa_f+\kappa_g)dx\llcorner B_R)}\nonumber\\
&\le W((\kappa_\lambda-\kappa_f)dx\llcorner B_R+\bar g,
(\kappa_\lambda-\kappa_f+\kappa_g)dx\llcorner B_R)
+W(\bar g,g').
\end{align}
According to Lemma \ref{lem:comp2}, we have for the first term
\begin{align}\label{as23}
W((\kappa_\lambda-\kappa_f)dx\llcorner B_R+\bar g,
&(\kappa_\lambda-\kappa_f+\kappa_g)dx\llcorner B_R)\le
\frac{1}{\kappa_\lambda-\kappa_f}\int_{B_R}|\nabla\bar\phi|^2,
\end{align}
where $\bar\phi$ is defined through
\begin{align*}
-\triangle\bar\phi=\kappa_g\;\;\mbox{in}\;B_R,\quad
\nu\cdot\nabla\bar\phi=\bar g\;\;\mbox{on}\;\partial B_R.
\end{align*}
We recall (\ref{ho08}),
which we may apply since $g'$ $\le\kappa_\lambda$ is supported in $\bar B_5-B_R$,
and which takes the form of
\begin{align}\label{as19}
\int_{\partial B_R}\bar g^2\le 2\cdot 5^{d-1}\kappa_\lambda\int||x|-R|dg'.
\end{align}
Combining this with (\ref{as24}), where $\phi$ is replaced by $\bar\phi$,
we obtain
\begin{align}\label{as16}
\int_{B_R}|\nabla\bar\phi|^2\le 2\cdot 5^{d-1}C_P\kappa_\lambda\int||x|-R|dg'.
\end{align}

\medskip

Before continuing with the r.~h.~s.~of (\ref{as16}),
we turn to the last term in (\ref{as14}). In view of (\ref{as15}), 
$\int\zeta(R\frac{x}{|x|},x)g'(dx)$ defines an admissible plan, so that
$W^2(\bar g,g')$ $\le\int|R\frac{x}{|x|}-x|^2g'(dx)$. Noting that
$|R\frac{x}{|x|}-x|$ $=||x|-R|$ and recalling that $g'$ is supported in $\bar B_5-B_R$,
we obtain
\begin{align*}
W^2(\bar g,g')\le 5\int||x|-R|dg'.
\end{align*}
In view of this and (\ref{as16}), we are lead to estimate $\int||x|-R|dg'$.
A simplification of the argument leading to (\ref{ho07}) gives
\begin{align}\label{as20}
\int_2^3\int||x|-R|dg' dR\le D.
\end{align}

\medskip

Since in view of (\ref{as10}) we have $\kappa_R$ $=\kappa_\lambda-\kappa_f+\kappa_g$, 
in order to obtain (\ref{as21}), it remains to show $\kappa_f^2+\kappa_g^2\lesssim D$.
Because of our assumption $D\ll 1$, this also deals with the pre-factors in (\ref{as11})
and (\ref{as23}).
By symmetry, we may restrict to $\kappa_g$. By definitions (\ref{as18}), (\ref{as17}),
and (\ref{as15}) we have $\kappa_g=\frac{1}{|B_R|}\bar g(\mathbb{R}^d)$, and thus by
Cauchy-Schwarz $\kappa_g^2\le\frac{|\partial B_R|}{|B_R|^2}\int_{\partial B_R}\bar g^2$.
Hence it remains to appeal to (\ref{as19}) and (\ref{as20}).


\subsection{A final approximation and proof of
Proposition \ref{prop:harmonicApproximation}}\label{PDEestimates}

Note that the two terms in (\ref{ao10}) involving $\nabla\phi(X(t))$
make only sense for general $\pi$ provided $\nabla\phi\in C^0(\bar B_R)$,
which however is not ensured by $\bar g-\bar f\in L^2(\partial B_R)$ in (\ref{ao92}).
Hence, a final -- however more conventional -- approximation argument is unavoidable:
We approximate $\bar g-\bar f$ by its mollification $(\bar g-\bar f)_r$
on a scale $r>0$, and denote by $\phi^r$ the corresponding solution of (\ref{ao92}).
With this replacement, (\ref{ao10}) assumes the form
\begin{align}
\lefteqn{\int_\Omega\int_\sigma^\tau|\dot X(t)-\nabla\phi^r(X(t))|^2dtd\pi}\label{as36}\\
&\le\int_{\Omega}|x-y|^2d\pi-\int_{B_R}|\nabla\phi^r|^2\label{as28}\\
&+2\int_{B_R}\phi^r\big(\frac{\mu(B_R)}{|B_R|}-\frac{\lambda(B_R)}{|B_R|}
-d(\mu-\lambda)\big)\label{as32}\\
&+2\int_{\partial B_R}\phi^r\big((\bar g-\bar f)_r-d(g-f)\big)\label{as30}\\
&+\int_\Omega\int_\sigma^\tau|\nabla\phi^r(X(t))|^2dtd\pi
-\int_{B_R}|\nabla\phi^r|^2.\label{as35}
\end{align}
Before controlling the difference in line (\ref{as35}), for which the mollification was made,
we address its effect on the lines
(\ref{as28}) (where it is cumbersome), (\ref{as32}) (where it is beneficial), 
and (\ref{as30}) (where it is both). 

\medskip

We now fix the radius $R$ such as to benefit from the Lemmas 
\ref{lem:approx}, \ref{lem:linfty}, and \ref{lem:localiseD}.
More precisely, taking the sum of the estimates (\ref{ao97}) divided by $o(E)+O(D)$,
(\ref{ao96}) divided by $E+D$, (\ref{m05}) by $o(E+D)$, (\ref{m06}) by $o(1)$, and 
(\ref{as21}) by $D$,\footnote{To pick out just the two terms of (\ref{ao97}) and (\ref{m06}), what we mean is the following: According to these to statements, for any given $\delta>0$ we have for sufficiently small $E+D$ that
\begin{align*}
\int_2^3dR \frac{1}{\delta E+D}W^2(f_R,\bar f_R)+\frac{1}{\delta}\pi(\Omega\cap\{\exists t\in[0,1]\;X(t)\in\partial B_R\})\lesssim 1.
\end{align*}
This means that there exists an $R\in[2,3]$ such that the integrand
is $\le 1$, which translates into
\begin{align*}
W^2(f_R,\bar f_R)\lesssim\delta E+D\quad\mbox{and}\quad
\pi(\Omega\cap\{\exists t\in[0,1]\;X(t)\in\partial B_R\})\lesssim\delta,
\end{align*}
which amounts to (\ref{cw08}) and (\ref{cw12bis}).}
we learn that there exists an
$R\in[2,3]$ and $\bar g,\bar f\in L^2(\partial B_R)$ with
\begin{align}
W^2(f,\bar f)+W^2(g,\bar g)&=o(E)+O(D),\label{cw08}\\
\int_{\partial B_R} \bar f^2+\bar g^2&=O(E+D),\label{cw06}\\
\int_{\Omega\cap\{\exists t\in[0,1]\; X(t)\in \partial B_R\}}|x-y|^2 d\pi&=o(E+D),\label{cw12}\\
\pi(\Omega\cap\{\exists t\in[0,1]\; X(t)\in \partial B_R\})&=o(1),\label{cw12bis}\\
W^2(\lambda\llcorner B_R,\frac{\lambda(B_R)}{|B_R|}dx\llcorner B_R)+(\frac{\lambda(B_R)}{|B_R|}-1)^2
&=O(D),\label{cw10}\\
W^2(\mu\llcorner B_R,\frac{\mu(B_R)}{|B_R|}dx\llcorner B_R)+(\frac{\mu(B_R)}{|B_R|}-1)^2
&=O(D).\label{cw11}
\end{align}

\medskip

We begin with line (\ref{as28}): 
It follows from a standard mollification argument that the
$\dot H^{-\frac{1}{2}}(\partial B_R)$-norm of $(\bar g-\bar f)_r-(\bar g-\bar f)$
is $O(r^\frac{1}{2})$ times the $L^2(\partial B_R)$-norm of $\bar g-\bar f$:
\begin{align}\label{as31}
\int_{\partial B_R}\zeta\big((\bar g-\bar f)_r-(\bar g-\bar f)\big)
\lesssim \big(r\int_{B_R}|\nabla\zeta|^2
\int_{\partial B_R}(\bar g-\bar f)^2\big)^\frac{1}{2}.
\end{align}
Hence from (\ref{ao48}) applied to (\ref{ao92}), which yields that
$\int_{B_R}|\nabla(\phi^{r}-\phi)|^2$ $=\int_{\partial B_R}(\phi^{r}-\phi)
((\bar g-\bar f)_r-(\bar g-\bar f))$, with the choice $\zeta = \phi^r-\phi$ in \eqref{as31}, we obtain
\begin{align}\label{as29}
\int_{B_R}|\nabla(\phi^{r}-\phi)|^2\lesssim r\int_{\partial B_R}(\bar g-\bar f)^2.
\end{align}
Combining this with (\ref{as24}) yields
\begin{align*}
\int_{B_R}|\nabla\phi|^2-\int_{B_R}|\nabla\phi^r|^2\leq 2 \int_{B_R}\nabla(\phi-\phi^r)\cdot \nabla\phi\lesssim 
r^\frac{1}{2}\int_{\partial B_R}(\bar g-\bar f)^2.
\end{align*}
Hence in view of (\ref{cw06}) we have
\begin{align}\label{cw07}
\int_{B_R}|\nabla\phi|^2-\int_{B_R}|\nabla\phi^r|^2\le r^\frac{1}{2}O(E+D).
\end{align}

\medskip

Equipped with (\ref{cw07}) and the more obvious
\begin{align}\label{cw09}
\int_{B_R}|\nabla\phi|^2=O(E+D),
\end{align}
we may now conclude the estimate of line (\ref{as28}).
We refer back to Subsection \ref{sec:harm}, namely
(\ref{cw05}) followed by Young's inequality. Next to an $o(E)$ coming from
the first factor in Young's inequality, this gives rise
to the square of the term in (\ref{ao94}) and the three terms in (\ref{ao95}). According to 
(\ref{cw10}) \& (\ref{cw11}) and (\ref{cw09}), the former contribution
is $O(D(E+D))$. According to (\ref{cw12}), once more (\ref{cw10}) \& (\ref{cw11}),
and (\ref{cw08}),
the latter three contributions in (\ref{ao95}) are $o(E+D)$, $O(D)$, and $o(E+D)$,
respectively. Hence in combination with (\ref{cw07}) we obtain
\begin{align*}
\lefteqn{\int_\Omega|x-y|^2d\pi-\int_{B_R}|\nabla\phi|^2}\nonumber\\
&\le o(E)+O(D(E+D))+o(E+D)+O(D)+r^\frac{1}{2}O(E+D),
\end{align*}
which is $o(E)+O(D)$ provided $r$ goes to zero as $E$ does.

\medskip

We now turn to the terms in the two lines (\ref{as32}) and (\ref{as30}).
Combining (\ref{as24}) with (\ref{as29}), we obtain by (\ref{cw06})
\begin{align}\label{as42}
\int_{B_R}|\nabla\phi^r|^2\lesssim\int_{\partial B_R}(\bar g-\bar f)^2.
\end{align}
Hence in conjunction with (\ref{as31}) we learn that 
\begin{align}\label{as33}
\mbox{in line (\ref{as30}) we may replace $(\bar g-\bar f)_r$ by $\bar g-\bar f$,}
\end{align}
once more at the expense of an error $r^\frac{1}{2}O(E+D)$. Now comes the
beneficial effect of mollification: By standard regularity theory for the 
Neumann-Poisson problem, $\sup_{\bar B_R}|\nabla\phi^r|$
is estimated by a sufficiently high norm of $(\bar g-\bar f)_r$, which
due to the mollification is controlled by $(\int_{\partial B_R}(\bar g-\bar f)^2)^\frac{1}{2}$.
A closer inspection shows that in line with scaling, this estimate assumes the form
\begin{align}\label{as34}
\sup_{\bar B_R}|\nabla\phi^{r}|^2\lesssim\frac{1}{r^{d-1}}\int_{\partial B_R}(\bar g-\bar f)^2
\stackrel{(\ref{cw06})}{=}\frac{1}{r^{d-1}}O(E+D).
\end{align}
Thus we may appeal to the following easy
consequence of the definition of $W^2$ (and Cauchy-Schwarz)
\begin{align}\label{as43}
\big|\int_{\partial B_R}\phi^r(\bar g-dg)\big|
\le\sup_{\bar B_R}|\nabla\phi^r|
\Big((\int_{\partial B_R}\bar g+g(\partial B_R))W^2(\bar g,g)\Big)^\frac{1}{2},
\end{align}
and a corresponding estimate for $f$.
Hence we learn from both (\ref{cw08}) and (\ref{cw06}) that the contribution from
(\ref{as33}) is $r^{-\frac{d-1}{2}}(o(E)+O(D))$.

\medskip

In conclusion, we obtain that the term in line (\ref{as30}) satisfies
\begin{align*}
\int_{\partial B_R}\phi^r((\bar g-\bar f)_r-d(g-f))
=\frac{1}{r^{\frac{d-1}{2}}}(o(E)+O(D))+r^\frac{1}{2}O(E+D),
\end{align*}
which still is $o(E)+O(D)$ provided $r$ only slowly goes to zero as $E+D$ does.
Appealing once more to an estimate similar to (\ref{as43}), but this time in conjunction with
(\ref{cw10}) and (\ref{cw11}), we learn that the term in line (\ref{as32}) is
\begin{align*}
\int_{B_R}\phi^r\big(\frac{\mu(B_R)}{|B_R|}-\frac{\lambda(B_R)}{|B_R|}-d(\mu-\lambda)\big)
=\frac{1}{r^{\frac{d-1}{2}}}O(D),
\end{align*} 
and thus is well-behaved, too.

\medskip

It remains to address the term in line \eqref{as35}; we claim that
\begin{align}\label{cw03}
\int_\Omega \int_\sigma^\tau |\nabla \phi^r(X(t))|^2d t d\pi&
-\int_{B_R} |\nabla \phi^r|^2\nonumber\\
&\le \frac 1 {r^{d-1}}o(E+D)+\frac 1 {r^d}(E+D)^\frac{3}{2}.
\end{align}
Indeed, setting $\zeta=|\nabla \phi^r|^2$, and appealing to (\ref{ao00bis}) and (\ref{ao66}),
we split the l.~h.~s.~of (\ref{cw03}) into the three differences
\begin{align}\label{as41}
\int_0^1\int_\Omega&\big(I(X(t)\in\bar B_R)\zeta(X(t))
-I(X(0)\in B_R)\zeta(X(0))\big)d\pi dt\nonumber\\
&+\left(\int_{B_R}\zeta d\lambda-\kappa_R\int_{B_R} \zeta\right)
+(\kappa_R-1)\int_{B_R} \zeta.
\end{align}
We turn to the first difference in \eqref{as41} and note
\begin{align*}
\lefteqn{I(X(t)\in \bar B_R)\zeta(X(t))-I(X(0)\in B_R)\zeta(X(0))}\\
&\leq I\big(\exists s\in[0,1]\; X(s)\in \partial B_R, X(t)\in\bar B_R\big)\zeta (X(t))\\
&+I\big(\forall s\in[0,1]\; X(s)\in B_R\big)\big(\zeta(X(t))-\zeta(X(0))\big),
\end{align*}
so that it is $\le$
\begin{align}\label{as40}
\sup_{\bar B_R}\zeta
\pi\big(\Omega\cap\{\exists s\in[0,1]\; X(s)\in \partial B_R\}\big)
+\sup_{\bar B_R}|\nabla \zeta|\int_\Omega|x-y|d\pi.
\end{align}
According to (\ref{cw12bis}), the second factor of the first term is $o(1)$.
The second factor of the second term is $\le(\pi(\Omega)\int_\Omega|x-y|^2d\pi)^\frac{1}{2}$, 
and thus $O(E^\frac{1}{2})$ by definitions (\ref{ao17}) and (\ref{ao45}).
Hence it remains to control the two first factors in (\ref{as40}).
Recalling the definition of $\zeta$ and \eqref{as34}, 
we have $\sup_{\bar B_R}\zeta$ $=r^{-(d-1)}O(E+D)$, so that the first term
in (\ref{as40}) is $r^{-(d-1)}o(E+D)$.
By the same argument that led to \eqref{as34}, we have
\begin{align}\label{cw14}
\sup_{\bar B_R}|\nabla^2 \phi^r|^2 \lesssim 
\frac{1}{r^{d+1}}\int_{\partial B_R}(\bar g-\bar f)^2
\stackrel{(\ref{cw06})}{=}\frac{1}{r^{d}}O(E+D),
\end{align}
which in combination with \eqref{as34} yields
\begin{align}\label{as44}
\sup_{\bar B_R} |\nabla \zeta|\lesssim \frac 1 {r^d} \int_{\partial B_R} (\bar g-\bar f)^2
=\frac{1}{r^d}O(E+D),
\end{align}
so that the second term in (\ref{as40}) is $r^{-d}O(E^\frac{1}{2}(E+D))$.

\medskip

We now address the second difference in (\ref{as41}). In analogy to (\ref{as43}), we
have that it is $\le$
\begin{align*}
\sup_{\bar B_R}|\nabla \zeta|\big(\lambda(B_R)+\kappa_R|B_R|) 
W^2( \lambda\llcorner B_R,\kappa_R dx\llcorner B_R)\big)^\frac{1}{2}.
\end{align*}
Due to (\ref{as21}) in Lemma \ref{lem:approx} and to \eqref{as44}
this is $r^{-d}O((E+D)D^\frac{1}{2})$.
Finally, the third term in \eqref{as41} is $\le(\sup\zeta)|\kappa_R-1|$
and thus $r^{-(d-1)}O((E+D)D^\frac{1}{2})$.

\medskip

The final task left is to estimate the l.~h.~s. of (\ref{ao89}) by the l.~h.~s. (\ref{as36}).
We start by approximating the argument of $\nabla\phi^r$:
\begin{align*}
\int_\Omega\int_\sigma^\tau|\nabla\phi^r(X(t))-I(x\in B_R)\nabla\phi^r(x)|^2dtd\pi.
\end{align*}
Using 
\begin{align*}
\lefteqn{I(t\in[\sigma,\tau])\big|\nabla\phi^r(X(t))-I(x\in B_R)\nabla\phi^r(x)\big)\big|^2}\\
&\stackrel{(\ref{ao66})}{=}\big|I(X(t)\in \bar B_R)\big(\nabla\phi^r(X(t))
-I(X(0)\in B_R)\nabla\phi^r(X(0))\big)\big|^2\\
&\leq I\big(\exists s\in[0,1]\; X(s)\in \partial B_R, X(t)\in\bar B_R\big)|\nabla\phi^r(X(t))|^2\\
&+I\big(\forall s\in[0,1]\; X(s)\in B_R\big)\big|\nabla\phi^r(X(t))-\nabla\phi^r(X(0))\big|^2,
\end{align*}
we obtain
\begin{align*}
\int_\Omega\int_\sigma^\tau&|\nabla\phi^r(X(t))-I(x\in B_R)\nabla\phi^r(x)|^2
dtd\pi\\
&\le\sup_{\bar B_R}|\nabla\phi^r|^2\pi\big(\Omega\cap\{\exists s\in[0,1]\; 
X(s)\in\partial B_R\}\big)\nonumber\\
&+\sup_{\bar B_R}|\nabla^2\phi^r|^2\int_\Omega|x-y|^2d\pi.
\end{align*}
By definition (\ref{ao17}) and (\ref{ao45}), as well as (\ref{cw12}), (\ref{as34}) and (\ref{cw14}) 
this yields
\begin{align}\label{cw15}
\int_\Omega\int_\sigma^\tau&|\nabla\phi^r(X(t))-I(x\in B_R)\nabla\phi^r(x)|^2
dtd\pi\nonumber\\
&\le \frac{1}{r^{d-1}}o(E+D)+\frac{1}{r^d}O((E+D)E).
\end{align}

\medskip

By the triangle inequality we obtain from (\ref{cw15})
\begin{align*}
\lefteqn{\int_\Omega(\tau-\sigma)|(y-x)-I(x\in B_R)\nabla\phi^r(x)|^2d\pi}\nonumber\\
&\le2\int_\Omega\int_\sigma^\tau|\dot X(t)-\nabla\phi^r(X(t))|^2dtd\pi\nonumber\\
&+\frac{1}{r^{d-1}}o(E+D)+\frac{1}{r^d}O((E+D)E).
\end{align*}
By definition (\ref{ao17}) of $\Omega$ and $R\ge 1$, we clearly have
\begin{align*}
\big((B_1\times\mathbb{R}^d)\cup(\mathbb{R}^d\times B_1)\big)\cap{\rm supp}\pi
\subset\Omega.
\end{align*}
Because of $R\ge 2$ and in our regime of short trajectories, cf.~(\ref{m07}), we have
\begin{align*}
\tau-\sigma=1,\;x\in B_R\quad
\mbox{for}\;(x,y)\in\big((B_1\times\mathbb{R}^d)\cup(\mathbb{R}^d\times B_1)\big)\cap{\rm supp}\pi.
\end{align*}
This implies
\begin{align*}
\lefteqn{\int_{(B_1\times\mathbb{R}^d)\cup(\mathbb{R}^d\times B_1)}
|(y-x)-\nabla\phi^r(x)|^2d\pi}\nonumber\\
&\le\int_\Omega(\tau-\sigma)|(y-x)-I(x\in B_R)\nabla\phi^r(x)|^2d\pi.
\end{align*}
The combination of the two last inequalities connects the l.~h.~s. of (\ref{ao89})
(with $\phi$ replaced by $\phi^r$) with the l.~h.~s. (\ref{as36}).
This concludes the proof of Proposition \ref{prop:harmonicApproximation}.


\end{document}